\documentclass[10pt, oneside, openright]{article}
\usepackage{epsfig}
\usepackage{graphicx, epstopdf}
\usepackage{graphpap}
\usepackage{amsmath, bm}
\usepackage{amsfonts}
\usepackage{amssymb}
\usepackage{makeidx}
\usepackage{amsmath,graphicx,url,times, amssymb, bm, epstopdf}
\usepackage{amssymb,amsmath,graphicx,times,url, mathtools}

\newcommand{\etab}{\mbox{\boldmath $\eta$}}
\newcommand{\psib}{\mbox{\boldmath $\psi$}}

\begin{document}

{\centering
{\Large Conservative Numerical Methods for Nonlinear String Dynamics: Non-planar Vibration}\\

\vspace{0.5in}
Technical Report, Acoustics and Audio Group, University of Edinburgh\\

\vspace{0.5in}

Stefan Bilbao\\

\vspace{0.5in}

May 9 2006\\}

\vspace{0.5in}

Nonlinear string vibration, in particular the case of nonplanar motion, has been an area of intense study for many years. Numerical simulation methods, essential for the comparison between measured data and theory, have received somewhat less attention. In this article, various numerical schemes for nonlinear nonplanar string dynamics are presented, with an emphasis on discrete conservation of energy and angular momentum. Simple numerical stability conditions may be arrived at, even under strongly nonlinear conditions, by employing these conservation properties. Full implementation details and various numerical examples are presented, and several topics, including a discussion of numerical loss models and spectral methods, are dealt with in the Appendix. \\

\newpage
\clearpage

\renewcommand{\baselinestretch}{1}
\begin{table}
\caption{\label{tab0} Nomenclature Table.} 
\begin{scriptsize}
\begin{tabular}{l|l}
Symbol & definition\\\hline
$A$ & string cross-sectional area\\
${\bf A}_{{\bf s}^{(\bullet)}}$ & update matrix, scheme ${\bf s}^{(\bullet)}$\\
${\mathcal A}_{{\bf S}}$, ${\mathcal A}_{{\bf K}}$, ${\mathcal A}_{{\bf s}^{(\bullet)}}$, ${\mathcal A}_{{\bf k}^{(\bullet)}}$ & angular momentum, model ${\bf S}$, ${\bf K}$, scheme ${\bf s}^{(\bullet)}$, ${\bf k}^{(\bullet)}$\\
${\bf B}_{{\bf s}^{(\bullet)}}$ & update matrix, scheme ${\bf s}^{(\bullet)}$\\
${\mathcal B}_{{\bf S},{\mathcal H}}$, ${\mathcal B}_{{\bf K},{\mathcal H}}$ & boundary term, energy conservation of system ${\bf S}$, ${\bf K}$\\
${\mathcal B}_{{\bf S},{\mathcal A}}$, ${\mathcal B}_{{\bf K},{\mathcal A}}$ & boundary term, angular momentum conservation of system ${\bf S}$, ${\bf K}$\\
${\mathcal B}_{{\bf s}^{(\bullet)},{\mathcal H}}$, ${\mathcal B}_{{\bf k}^{(\bullet)},{\mathcal H}}$ & boundary term, energy conservation of system ${\bf s}^{(\bullet)}$, ${\bf k}^{(\bullet)}$\\
${\mathcal B}_{{\bf s}^{(\bullet)},{\mathcal A}}$, ${\mathcal B}_{{\bf k}^{(\bullet)},{\mathcal A}}$ & boundary term, angular momentum conservation of system ${\bf s}^{(\bullet)}$, ${\bf k}^{(\bullet)}$\\
${\bf D}$ & differentiation matrix\\
${\mathcal D}$ & discrete spatial domain\\
$e_{\bullet}$ & shift operator\\
$E$ & Young's modulus\\
$h_{t}$ & time step\\
$h_{x}$ & grid spacing\\
${\mathcal G}$, ${\mathcal G}_{{\bf k}^{(\bullet)}}$ & nonlinear scaling, system ${\bf K}$, scheme ${\bf k}^{(\bullet)}$\\
${\mathcal H}_{{\bf S}}$, ${\mathcal H}_{{\bf K}}$, ${\mathcal H}_{{\bf s}^{(\bullet)}}$, ${\mathcal H}_{{\bf k}^{(\bullet)}}$  & total energy, model ${\bf S}$, ${\bf K}$, scheme ${\bf s}^{(\bullet)}$, ${\bf k}^{(\bullet)}$\\
$L$ & string length\\
$n$ & time index \\
$p$ & $=\xi_{x}$ until Section \ref{sksec}, $=\delta_{x-}\xi$ thereafter\\
${\bf P}$ & auxiliary scheme matrix \\
${\bf q}$ & $=\etab_{x}$ until Section \ref{sksec}, $=\delta_{x-}\etab$ thereafter\\
${\bf Q}_{(1)},{\bf Q}_{(2)}$ & auxiliary scheme matrices \\
$t$ & time variable (non-dimensionalized from Section \ref{ndsec} onwards)\\
\pagebreak
$T_{0}$ & nominal string tension\\
${\mathcal T}_{{\bf S}}$, ${\mathcal T}_{{\bf K}}$, ${\mathcal T}_{{\bf s}^{(\bullet)}}$, ${\mathcal T}_{{\bf k}^{(\bullet)}}$ & kinetic energy, model ${\bf S}$, ${\bf K}$, scheme ${\bf s}^{(\bullet)}$, ${\bf k}^{(\bullet)}$\\
${\bf u}$ & vector of discrete longitudinal displacements\\
${\bf v}_{(1)},{\bf v}_{(2)}$ & vectors of discrete transverse displacements\\
${\mathcal V}_{{\bf S}}$, ${\mathcal V}_{{\bf K}}$, ${\mathcal V}_{{\bf s}^{(\bullet)}}$, ${\mathcal V}_{{\bf k}^{(\bullet)}}$ & potential energy, model ${\bf S}$, ${\bf K}$, scheme ${\bf s}^{(\bullet)}$, ${\bf k}^{(\bullet)}$\\
${\bf w}$ & discrete state update vector \\
$x$ & spatial independent variable (non-dimensionalized from Section \ref{ndsec} onwards)\\
$\alpha$ & dimensionless string parameter $( =T_{0}/EA$)\\
$\beta$ & = $\frac{\alpha-1}{2}$\\
$\delta_{\bullet}$ & difference operator\\
$\eta_{(1)},\eta_{(2)}$ & transverse displacements\\
$\etab$ & vector transverse displacement (non-dimensionalized from Section \ref{ndsec} onwards)\\
$\hat{\etab}$ & Fourier expansion coefficients of vector transverse displacement \\
$\lambda$ & Courant number $(=h_{t}/h_{x})$\\
$\mu_{\bullet}$ & discrete averaging operator \\
$\nu$ & free parameter, implicit scheme \\
$\xi$ & longitudinal displacement (non-dimensionalized from Section \ref{ndsec} onwards)\\
$\rho$ & linear mass density\\
$\sigma_{\xi}$, $\sigma_{{\scriptsize\etab}}$ & loss parameter (longitudinal, transverse)\\
$\tau$ & free parameter, implicit scheme \\
$\phi$ & $=p+\frac{1-\alpha}{2}{\bf q}^{T}{\bf q}$\\
$\psib$ & $=\!\left(\alpha+\frac{1-\alpha}{2}\left({\bf q}^{T}{\bf q}+2p\right)\right){\bf q}$\\
 \end{tabular}
\end{scriptsize}
\end{table}

\newpage\clearpage

\section{Introduction}
\label{introsec}

The dynamics of strings under nonlinear conditions has been under study at least since the early work of Kirchhoff \cite{Kirchhoff} and then Carrier \cite{Carrier}; later, Anand \cite{Anand} and  Narasimha \cite{Narasimha} extended this work considerably.  This easily defined, but nonetheless extremely complex system exhibits a wide variety of behaviour characteristic of nonlinear systems; of particular interest is the ``whirling" phenomenon, which is peculiar to the case of non-planar motion \cite{Johnson}, \cite{Rubin}. Nonlinear string dynamics in three dimensions serves as an excellent test problem not only for a variety of analysis techniques, but also for the construction of numerical methods, which are the focus of this article.

The study of numerical methods which inherit discrete conservation laws from a continuous model system has been ongoing for some time \cite{Sanz82}, \cite{Greenspan}, and relates to early work on the so-called energy-method \cite{Richtmyer}. Most often, the systems under study are of a general form, e.g., the nonlinear Klein-Gordon equation, which was approached by Vu-quoc and Li \cite{vuquoc93} and Li and Vu-quoc \cite{vuquoc}, and various systems including a single polarization transverse-only model of string dynamics which were discussed by Furihata \cite{Furihata}. In all these cases, the form of the nonlinearity is left unspecified. In various commonly-encountered models of string dynamics, the nonlinearity is often simplified using a series approximation. Such a simplification allows obvious benefits in the analysis of the string, particularly when the nonlinearity is approximated using quadratic or cubic terms \cite{Nayfeh}, as it often is. At the same time, the possibility of exploiting various algebraic symmetries in a nonlinear difference scheme also appears; as will be shown in this article, there are many distinct ways of designing numerical schemes for the same nonlinear string system, which vary considerably in terms of their conservation properties, ease of use, and, most importantly, their stability properties. Conservative difference schemes for the Kirchhoff-Carrier string, and for a more general string undergoing planar motion have been discussed by this author in \cite{Bilbaoacustica1} and \cite{Bilbaojasa05}, respectively. 

It is perhaps worth mentioning here that one extremely interesting recent application of numerical simulation techniques for nonlinear mechanical systems such as the string is in the area of musical sound synthesis. Such physical modeling synthesis, as it is often called, has been in existence for some time now, and for some systems, real-time performance is now possible using personal computers. There has been some important work on using standard numerical methods such as finite difference schemes \cite{Ruiz}, \cite{Bacon}, \cite{Chaigne92}, but the dominant techniques have been based around efficient structures with their roots in digital filter design; probably the best known are digital waveguides \cite{josbook}. Recently, there have been some attempts at sound synthesis based on nonlinear string vibration, in order to model effects such as pitch glides under high amplitude plucking conditions \cite{Vesa99}, as well as the phantom partial phenomenon which occurs in piano strings \cite{Bankdafx04}. Though this article is intended for a general audience, some commentary on this topic will appear at various points throughout this article. 

In Section \ref{modelsec}, a general model of nonlinear string dynamics is presented, followed by two simpler forms, one employing the series approximations to the nonlinearity mentioned above (system {\bf S}), and a further simplified form of the Kirchhoff-Carrier variety (system {\bf K}). The section concludes with a brief presentation of the energy and angular momentum conservation properties of these models, and in particular, the bounds on the growth of the solution which result from the former property. Section \ref{fdsec} is a short recap of the properties of finite difference operators and inner product spaces, with a view toward applications in the construction of conservative schemes. In Section \ref{sksec}, various difference approximations to systems {\bf S} and {\bf K} are presented, followed by an analysis of their discrete energy and angular momentum conservation properties. The discrete conserved energy for each of the schemes is further examined, first to determine conditions under which, when it indeed exists, it is positive, and if so, what bounds may be placed on the size of the solution. The section concludes with a brief look at the schemes in the forms in which they will be implemented. Numerical results are presented in Section \ref{numsec}, with a special emphasis on phenomena which are inherent to motion in three dimensions, and in particular the so-called ``whilrling" behaviour \cite{Johnson}, \cite{Rubin}. Several other topics are briefly addressed in the Appendices, namely the generalization to the case of linear damping in Appendix \ref{losssec}, a loosening of the stability condition on the time step in Appendix \ref{impsec}, and finally, in Appendix \ref{specsec}, a brief look at a spectral-type method for the integration of system {\bf K}.

\section{Nonlinear String Models}
\label{modelsec}

\begin{figure}[ht]
\centerline{\includegraphics[scale=0.75,clip,trim=14mm 0mm 58mm 166mm]{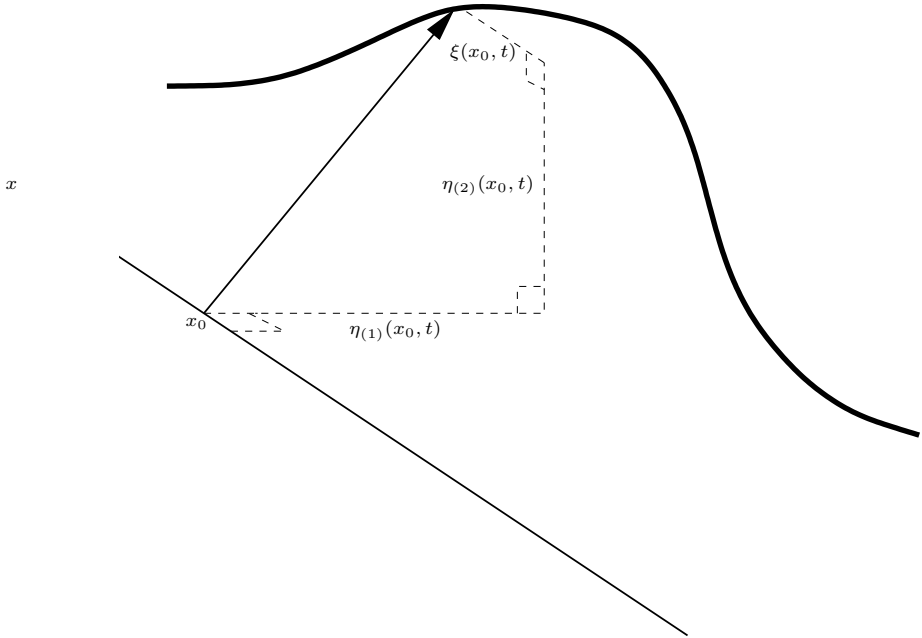}\put(-182,220){{\scriptsize $\xi(x_{0},t)$}}
\put(-220,115){{\scriptsize $\eta_{(1)}(x_{0},t)$}}\put(-350,170){{\scriptsize $x$}}\put(-282,118){{\scriptsize $x_{0}$}}\put(-185,170){{\scriptsize $\eta_{(2)}(x_{0},t)$}}}
\caption{Geometry of three-dimensional string displacement, illustrating displacement of string from rest point $(x_{0},0,0)$ to $(x+\xi(x_{0},y), \eta_{(1)}(x_{0},t),\eta_{(2)}(x_{0},t))$. }  
\label{fig0}
\end{figure}

A general model of nonlinear string dynamics, discussed by many authors, and summarized succinctly by Morse and Ingard \cite{Morse} and which can be related to the geometrically-exact theory of beams \cite{simo87}, is given by the following set of equations:
\begin{subequations}
\label{Morsed}
\begin{eqnarray}
\rho\xi_{tt} \!&=& \!EA\xi_{xx}-(EA-T_{0})\!\!\left(\!\frac{1+\xi_{x}}{\sqrt{\left(1+\xi_{x}\right)^2+\etab_{x}^{T}\mbox{{\boldmath $\eta$}}_{x}}}\!\right)_{x}\\
\rho\etab_{tt}\! &=& \!EA\etab_{xx}-(EA-T_{0})\!\!\left(\!\frac{\etab_{x}}{\sqrt{\left(1+\xi_{x}\right)^2+\etab_{x}^{T}\etab_{x}}}\!\right)_{x}
\end{eqnarray}
\end{subequations}
Here, $\xi(x,t)$ and the two-element column vector $\etab(x,t) =  [\eta_{(1)},\eta_{(2)}]^{T}$ describe, respectively, the longitudinal and transverse deviation of a point on the string as a function of time $t\geq 0$ and distance along the string $x\in [0,L]$. (The superscript $^T$ indicates a vector transpose.)  Such a point, located at Cartesian coordinates $(x,0,0)$ when the string is at rest, will have dynamic coordinates $(x+\xi, \eta_{(1)}, \eta_{(2)})$. See Figure \ref{fig0}. $E$, $A$, $\rho$ and $T_{0}$ are Young's modulus, cross-sectional area, linear mass density, and nominal tension for the string, all assumed constant here. Subscripts indicate differentiation with respect to the named independent variable. System \eqref{Morsed} is by no means the most general model of string dynamics; higher order effects may be modelled as well, as per the work of Narasimha \cite{Narasimha} and Kurmyshev \cite{Kurmyshev}; such improved models may fall outside the range of the techniques presented here.  


System \eqref{Morsed} requires the specification of the initial conditions $\xi (x,0)$, $\xi_{t}(x,0)$, $ \etab(x,0)$, and $\etab_{t}(x,0)$, in order that the solution exist and be uniquely defined. A discussion of boundary conditions is postponed briefly until Section \ref{bcsec}. 

\subsection{A Nondimensionalized Form}
\label{ndsec}
System \eqref{Morsed} may be nondimensionalized by introducing the variables
\begin{equation}
x' = x/L\qquad \xi = \xi/L\qquad \etab = \etab/L \qquad t' = t\sqrt{\frac{EA}{\rho L^2}}
\end{equation}
which, when substituted in \eqref{Morsed} and primes removed, leads to
\begin{subequations}
\label{Morse}
\begin{eqnarray}
\xi_{tt} &=& \xi_{xx}-(1-\alpha)\left(\frac{1+\xi_{x}}{\sqrt{\left(1+\xi_{x}\right)^2+\etab_{x}^{T}\mbox{{\boldmath $\eta$}}_{x}}}\right)_{x}\\
\etab_{tt} &=& \etab_{xx}-(1-\alpha)\left(\frac{\etab_{x}}{\sqrt{\left(1+\xi_{x}\right)^2+\etab_{x}^{T}\etab_{x}}}\right)_{x}
\end{eqnarray}
\end{subequations}
which depends on a single parameter $\alpha = T_{0}/EA$, and which is defined over $x\in[0,1]$ 

\subsection{Approximate Forms}

There are various levels of approximation to system \eqref{Morse}; the most direct forms result from Taylor series approximations to the nonlinearity (i.e., the term in large parentheses). An approximation to first order uncouples the longitudinal and transverse motion (and the two transverse polarizations from one another), yielding linear wave equations, with wave speeds 1 (longitudinal) and $\sqrt{\alpha}$ (both transverse polarizations). An approximation to second order is sometimes employed \cite{Bankdafx04}, \cite{Banksmac03}, but most common in the study of nonlinear string vibration is an approximation to third order, following Anand \cite{Anand} and Morse \cite{Morse} in making use of the observation that $\xi = O(|\etab|^{2})$, which is true for metallic strings: 

{\bf System S}

\begin{subequations}
\label{Morse3}
\begin{eqnarray}
\label{Morsesys31}
\xi_{tt} \!&=& \phi_{x}\qquad\mbox{{\rm with}}\qquad\phi=p+\frac{1-\alpha}{2}{\bf q}^{T}{\bf q}\\
\label{Morsesys32}
\etab_{tt}\! &=& \psib_{x}\qquad\mbox{{\rm with}}\qquad\psib=\!\left(\alpha+\frac{1-\alpha}{2}\left({\bf q}^{T}{\bf q}+2p\right)\right){\bf q}
\end{eqnarray}
\end{subequations}
For notational simplicity, the symbols $p$ and ${\bf q}$ have been introduced; they are defined as
\begin{equation}
\label{pqcdef}
p = \xi_{x}\qquad {\bf q} = \etab_{x}
\end{equation}
Numerical methods for system ${\bf S}$, in its planar form, have been applied to the problem of piano string vibration at high amplitudes \cite{Bankdafx04}, in order to generate perceptually important ``phantom partials," \cite{Conklin} which result from coupling between longitudinal and transverse motion. 
 
Under certain conditions, namely that (1) $1\gg \alpha$, (2) the term $\xi_{tt}$ in Eq. \eqref{Morsesys31} may be neglected, and (3) the longitudinal displacement $\xi(x,t)$ is zero at $x=0$ and $x=1$, system {\bf S} may be reduced to a single equation in the vector transverse displacement \etab\, alone:

{\bf System K}

\begin{equation}
\label{kc}
\etab_{tt} = \alpha{\mathcal G}{\bf q}_{x}
\end{equation}
where ${\mathcal G}$ is defined by
\begin{equation}
\label{gdef}
{\mathcal G} = 1+\frac{1}{2\alpha}\int_{0}^{1}{\bf q}^{T}{\bf q}dx
\end{equation}

System {\bf K} above, often referred to as being of Kirchhoff-Carrier type \cite{Kirchhoff}, \cite{Carrier}, is far simpler to deal with that system {\bf S}, both analytically and numerically, for the simple reason that the nonlinearity, characterized by ${\mathcal G}$, is averaged over the string length, and does not have any spatial dependence (i.e., it is a scalar quantity). Eq. \eqref{kc} can be viewed, in a crude sense, as a wave equation with a wave speed which is dependent on variations in total string length. This type of system has been employed, in the context of digital waveguides \cite{josbook} in order to produce synthetic sound for strings under high-amplitude plucking conditions \cite{Vesa99}. This type of nonlinearity is often referred to as ``tension modulation" in the sound synthesis community. 

\subsection{Conserved Quantities}

Given systems {\bf S} and {\bf K}, which will serve as the models to be solved numerically in this article, it is worthwhile to spend some time examining their conservation properties. In the present case of continuously variable systems, this is quite straightforward. 

\subsubsection{Energy Conservation}

For system {\bf S}, multiplying Eq. \eqref{Morsesys31} by $\xi_{t}$ and left-multiplying Eq. \eqref{Morsesys32} by \etab$^{T}_{t}$ and then integrating over the interval $[0,1]$ gives
\begin{equation*}
\int_{0}^{1}\xi_{t}\xi_{tt} -\xi_{t}\phi_{x}dx = 0\qquad\qquad\int_{0}^{1}\etab_{t}^{T}\etab_{tt}-\etab_{t}^{T}\psib_{x}dx=0
\end{equation*}

Integrating by parts for the latter term under the integral in each equation and adding the results gives, employing definitions Eqs. \eqref{pqcdef},
\begin{eqnarray}
\label{Segy}
\int_{0}^{1}\xi_{t}\xi_{tt} + \etab_{t}^{T}\etab_{tt}+p_{t}\phi + {\bf q}^{T}_{t}\psib dx={\mathcal B}_{{\bf S},{\mathcal H}}\notag
\end{eqnarray}
where the boundary term ${\mathcal B}_{{\bf S},{\mathcal H}}$ is given by
\begin{equation}
\label{bshdef}
{\mathcal B}_{{\bf S},{\mathcal H}} = \left(\xi_{t}\phi+\etab^{T}_{t}\psib\right)\Big|_{0}^{1}
\end{equation}
Eq. \eqref{Segy} can be written as 
\begin{equation*}
\frac{d}{dt}{\mathcal H}_{{\bf S}} = {\mathcal B}_{{\bf S},{\mathcal H}}\qquad\mbox{{\rm for}}\qquad {\mathcal H}_{{\bf S}} = {\mathcal T}_{{\bf S}}+{\mathcal V}_{{\bf S}}
\end{equation*}
where $d/dt$ represents a total derivative with respect to time $t$, and the scalar quantities ${\mathcal T}_{{\bf S}}$ and ${\mathcal V}_{{\bf S}}$ are given by
\begin{eqnarray*}
{\mathcal T}_{{\bf S}} &=& \int_{0}^{1}\frac{1}{2}\xi_{t}^{2}+\frac{1}{2}\etab^{T}_{t}\etab_{t}dx\\
{\mathcal V}_{{\bf S}} &=& \int_{0}^{1}\frac{\alpha}{2}p^{2}+\frac{\alpha}{2}{\bf q}^{T}{\bf q}+\frac{1-\alpha}{2}\left(p+\frac{1}{2}{\bf q}^{T}{\bf q}\right)^{2}dx
\end{eqnarray*}
Clearly, ${\mathcal T}_{{\bf S}}$ and ${\mathcal V}_{{\bf S}}$ represent the kinetic and potential energy, respectively, of the string as described by system {\bf S}, and ${\mathcal H}_{{\bf S}}$ the total energy, whose rate of change is dependent only on the boundary term ${\mathcal B}_{{\bf S},{\mathcal H}}$. 

For system {\bf K}, a similar property may be derived, i.e., 
\begin{equation*}
\frac{d}{dt}{\mathcal H}_{{\bf K}} = {\mathcal B}_{{\bf K},{\mathcal H}}\qquad\mbox{{\rm for}}\qquad {\mathcal H}_{{\bf K}} = {\mathcal T}_{{\bf K}}+{\mathcal V}_{{\bf K}}
\end{equation*}
with 
\begin{eqnarray*}
{\mathcal T}_{{\bf K}} &=& \int_{0}^{1}\frac{1}{2}\etab^{T}_{t}\etab_{t}dx\\
{\mathcal V}_{{\bf K}} &=& \frac{\alpha}{2}\left(\int_{0}^{1}{\bf q}^{T}{\bf q}dx\right)\left(1+\frac{1}{4\alpha}\int_{0}^{1}{\bf q}^{T}{\bf q}dx\right)\\
{\mathcal B}_{{\bf K},{\mathcal H}} &=& \alpha{\mathcal G}\etab^{T}_{t}{\bf q}\Big|_{0}^{1}
\end{eqnarray*}
\subsubsection{Angular Momentum Conservation}

In order to examine the conservation of angular momentum, it is useful to define, for any two-vector ${\bf f} = [f_{(1)}, f_{(2)}]^{T}$, the operation $\tilde{}$ by $\tilde{{\bf f}} = [-f_{(2)}, f_{(1)}]^{T}$. It should be clear that for any such vector ${\bf f}$, it must be true that 
\begin{equation}
\label{tildeidc}
\tilde{{\bf f}}^{T}{\bf f} = 0
\end{equation}

Considering only the angular momentum of system {\bf S} about the string axis (i.e., the $x$-component), left-multiply Eq. \eqref{Morsesys32} by $\tilde{\etab}^{T}$ and integrate over the interval $[0,1]$ to get
\begin{equation*}
\int_{0}^{1}\tilde{\etab}^{T}\etab_{tt}-\tilde{\etab}^{T}\psib_{x}dx=0
\end{equation*}
Expanding the first term and integrating the second by parts gives
\begin{equation}
\label{Sam}
\int_{0}^{1}\left((\tilde{\etab}^{T}\etab_{t})_{t}-\tilde{\etab}^{T}_{t}\etab_{t}\right)+\tilde{{\bf q}}^{T}\psib dx={\mathcal B}_{{\bf S},{\mathcal A}}
\end{equation}
where
\begin{equation*}
{\mathcal B}_{{\bf S},{\mathcal A}} = \tilde{\etab}^{T}\psib \Big|_{0}^{1}
\end{equation*}
Finally, by applying identity \eqref{tildeidc} to the second and third terms under the integral above (note that \psib\, is proportional to ${\bf q}$), Eq. \eqref{Sam} can be reduced to 
\begin{equation*}
\frac{d}{dt}{\mathcal A}_{{\bf S}} = {\mathcal B}_{{\bf S},{\mathcal A}}
\end{equation*}
with 
\begin{equation*}
{\mathcal A}_{{\bf S}} = \int_{0}^{1}\tilde{\etab}^{T}\etab_{t}dx
\end{equation*}
which is the total angular momentum of system {\bf S}, in the $x$ direction. 

For system {\bf K}, the analysis is similar, and gives 
\begin{equation*}
\frac{d}{dt}{\mathcal A}_{{\bf K}} = {\mathcal B}_{{\bf K},{\mathcal A}}
\end{equation*}
with 
\begin{eqnarray*}
{\mathcal A}_{{\bf K}} &=& \int_{0}^{1}\tilde{\etab}^{T}\etab_{t}dx\\
{\mathcal B}_{{\bf K},{\mathcal A}} &=& \alpha{\mathcal G}\tilde{\etab}^{T}{\bf q}\Big|_{0}^{1}
\end{eqnarray*}
${\mathcal A}_{{\bf K}}$ is the total angular momentum of system {\bf K}. 

\subsection{Boundary Conditions}
\label{bcsec}
In the interest of simplifying the analysis somewhat, and of rendering systems {\bf S} and {\bf K} conservative, it is useful to specify several boundary conditions. Fixed conditions, at an end of the string, are defined by
\begin{equation}
\label{bcdef1}
\xi = 0\qquad\mbox{{\rm and}}\qquad \etab= {\bf 0} \\
\end{equation}
and free conditions by
\begin{equation}
\label{bcdef2}
p = 0\qquad\mbox{{\rm and}}\qquad {\bf q}={\bf 0} \\
\end{equation}
Note that the free conditions above also imply that $\phi$ and \psib\, vanish. 

If either of these conditions holds at each of $x=0$ and $x=1$, then both ${\mathcal B}_{{\bf S},{\mathcal H}}$ and ${\mathcal B}_{{\bf S},{\mathcal A}}$ vanish, and as a result
\begin{equation}
\label{Segycons}
{\mathcal H}_{{\bf S}}(t) = {\mathcal H}_{{\bf S}}(0)\qquad {\mathcal A}_{{\bf S}}(t) = {\mathcal A}_{{\bf S}}(0)
\end{equation}
In other words, the energy and angular momentum of system {\bf S} remain constant for all $t$ and equal to their initial values. Many other conditions, including mixtures of the above conditions, obviously lead to conservative behaviour as well. 

For system {\bf K}, of the two conditions given above, only the conditions \eqref{bcdef1} above is allowable (recall that $\xi=0$ at either end of the string is a starting point in the derivation of system {\bf K}), and again leads, when enforced at both ends of the string, to   
\begin{equation}
\label{Kegycons}
{\mathcal H}_{{\bf K}}(t) = {\mathcal H}_{{\bf K}}(0)\qquad {\mathcal A}_{{\bf K}}(t) = {\mathcal A}_{{\bf K}}(0)
\end{equation}
\subsection{Bounds on Solution Size}
\label{boundsec}
In the interest of simplifying notation, it is useful to introduce the spatial $L_{2}$ inner product and norm of two column vector functions of $x\in[0,1]$ and perhaps $t$, ${\bf f}$ and ${\bf g}$, containing the same number of elements. These are given by
\begin{equation*}
\langle {\bf f},{\bf g}\rangle = \int_{0}^{1}{\bf f}^{T}{\bf g} dx\qquad \|{\bf f}\| = \langle {\bf f},{\bf f}\rangle^{1/2}
\end{equation*}
Obviously, such norms and inner products remain functions of the time variable $t$; when necessary, this will be indicated, e.g., $\|{\bf f}\|(t)$. 

Returning to the expressions for kinetic and potential energy, and angular momentum of system {\bf S}, it is then possible to rewrite them as
\begin{eqnarray*}
{\mathcal T}_{{\bf S}} &=& \frac{1}{2}\|\xi_{t}\|^{2}+\frac{1}{2}\|\etab_{t}\|^2\\
{\mathcal V}_{{\bf S}} &=& \frac{\alpha}{2}\|p\|^{2}+\frac{\alpha}{2}\|{\bf q}\|^2+\frac{1-\alpha}{2}\|p+\frac{1}{2}{\bf q}^{T}{\bf q}\|^{2}\\
{\mathcal A}_{{\bf S}} &=& \langle \tilde{\etab}, \etab_{t}\rangle
\end{eqnarray*}
${\mathcal T}_{{\bf S}}$ is non-negative, and ${\mathcal V}_{{\bf S}}$ is as well, under the condition
\begin{equation}
\label{EATcond}
\alpha\leq 1
\end{equation}
which is the case for moderately stiff strings, and which will be assumed henceforth in this article. 

Similarly, for system {\bf K}, the kinetic and potential energies and angular momentum may be written as
\begin{eqnarray*}
{\mathcal T}_{{\bf K}} &=& \frac{1}{2}\|\etab_{t}\|^2\\
{\mathcal V}_{{\bf K}} &=& \frac{\alpha}{2}\|{\bf q}\|^2\left(1+\frac{B}{2}\|{\bf q}\|^2\right)\\
{\mathcal A}_{{\bf K}} &=& \langle \tilde{\etab}, \etab_{t}\rangle
\end{eqnarray*}
\subsubsection{General Bounds}

For system {\bf S}, conservative boundary conditions lead to Eq. \eqref{Segycons}, and it is clearly true then that ${\mathcal T}_{{\bf S}}(t)\leq {\mathcal H}_{{\bf S}}(t)={\mathcal H}_{{\bf S}}(0)$, further implying that
\begin{equation}
\label{xietabound}
\|\xi_{t}\|(t)\leq\sqrt{2{\mathcal H}_{{\bf S}}(0)}\qquad\|\etab_{t}\|(t)\leq\sqrt{2{\mathcal H}_{{\bf S}}(0)}
\end{equation}

It is simple to translate the above bounds on the time derivatives of the displacements to bounds on the displacements themselves. For the longitudinal displacement, for example, one may write, employing the Cauchy-Schwartz inequality \cite{Horn} and the first of the bounds \eqref{xietabound},
\begin{equation*}
2\|\xi\|\frac{d}{dt}\|\xi\| = \frac{d}{dt}\|\xi\|^2 = 2\langle \xi, \xi_{t}\rangle \leq 2\|\xi\|\|\xi_{t}\|\leq 2\|\xi\|\sqrt{2{\mathcal H}_{{\bf S}}(0)}
\end{equation*}
implying that
\begin{equation*}
\frac{d}{dt}\|\xi\|\leq\sqrt{2{\mathcal H}_{{\bf S}}}\Longrightarrow\|\xi\|(t) \leq \|\xi\|(0)+ \sqrt{2{\mathcal H}_{{\bf S}}(0)}t
\end{equation*}
In other words, growth of the $L_2$ norm of the longitudinal displacement is at most linear. An identical bound follows for the transverse displacement, i.e., 
\begin{equation*}
\|\etab\|(t) \leq \|\etab\|(0)+ \sqrt{2{\mathcal H}_{{\bf S}}(0)}t
\end{equation*}
For system {\bf K}, a similar bound on the transverse displacement holds, under conservative boundary conditions:
\begin{equation*}
\|\etab\|(t) \leq \|\etab\|(0)+ \sqrt{2{\mathcal H}_{{\bf K}}(0)}t
\end{equation*}

\subsubsection{Bounds under Fixed Conditions}

For a given type of motion (i.e., longitudinal or transverse), if at least one of the ends is fixed, then better bounds are possible. For instance, consider condition \eqref{bcdef1} applied at $x=0$, paying special attention to the first condition, i.e., $\xi(0,t)=0$. The following bound is immediate:
\begin{equation*}
\xi = \int_{0}^{x}p dx \leq \left(\int_{0}^{x}p^2 dx\right)^{1/2}\left(\int_{0}^{x}1 dx\right)^{1/2}\leq\|p\|
\end{equation*}
which implies, furthermore, that
\begin{equation*}
\xi^2\leq \|p\|^2  \Longrightarrow \|\xi\|(t)\leq \|p\|(t)
\end{equation*}
(If $\xi = 0$ at both $x=0$ and $x=1$, then the above bound may be tightened to $\|\xi\|\leq \|p\|/2$.) Finally, using the fact that ${\mathcal V}_{{\bf S}}(t)\leq {\mathcal H}_{{\bf S}}(t)={\mathcal H}_{{\bf S}}(0)$, one may then conclude that
\begin{equation*}
\|p\|(t)\leq\sqrt{\frac{2{\mathcal H}_{{\bf S}}(0)}{\alpha}}\Longrightarrow \|\xi\|(t)\leq \sqrt{\frac{2{\mathcal H}_{{\bf S}}(0)}{\alpha}}
\end{equation*}
Similar bounds may be found for \etab, again under fixed transverse conditions at at least one end of the string, for both system {\bf S} and system {\bf K}. 

\section{Grid Functions and Finite Difference Operators}
\label{fdsec}
In this section, a short review of grid functions as employed in finite difference schemes and the related difference operators is provided, with an eye toward applications in the study of schemes for conservative systems. Due to the vector nature of the differential equations to be studied here, definitions are framed here in terms of vector-valued grid functions, but it should be kept in mind that in most cases, such definitions reduce simply to the scalar case. 

\subsection{Grid Functions}

A grid function ${\bf f}_{i}^{n}$ is defined as a column vector taking on values at the collection of points indexed by integers $i$ and $n\geq 0$; $({\bf f}_{i}^{n})^{T}$ is its transpose. It is intended as an approximation to a continuously-variable function ${\bf f}(x,t)$ at the location $x=ih_{x}$, $t=nh_{t}$, where $h_{x}$ and $h_{t}$ are the grid spacing and time step, respectively. In order that the connection with the underlying model problem be maintained, in this article, a grid function will be described using the same variable name as the continuous function it approximates. In addition, if a grid function is presented without one or both of its indices, it is assumed to have general indices $i$ and $n$. 

As for the continuous case, for any column two-vector grid function ${\bf f} = [f_{(1)},f_{(2)}]$, the grid function $\tilde{{\bf f}}$ is defined by $\tilde{{\bf f}} = [-f_{(2)},f_{(1)}]$. It is then always true that, for any two two-vector grid functions ${\bf f}$ and ${\bf g}$, 
\begin{equation}
\label{tildeflipid}
\tilde{{\bf f}}^{T}{\bf g} = -{\bf f}^{T}\tilde{{\bf g}}
\end{equation}
and, furthermore, 
\begin{equation}
\label{tildeid}
\tilde{{\bf f}}^{T}{\bf f} = 0
\end{equation}  
where here, ``0" is interpreted as a grid function taking on the value zero for all $i$ and $n$. 

The important parameter $\lambda$, defined by
\begin{equation}
\label{Courantdef}
\lambda = h_{t}/h_{x}
\end{equation}
is crucial in that numerical stability conditions are framed in terms of it, as per the linear case \cite{Strikwerda}. 

\subsection{Temporal Operators}
\label{timesec}
The basic temporal operators are the unit forward and backwards shifts, defined in terms of their effect on a grid function ${\bf f}$ by
\begin{equation*}
e_{t+}{\bf f}_{i}^{n} = {\bf f}_{i}^{n+1}\qquad e_{t-}{\bf f}_{i}^{n} = {\bf f}_{i}^{n-1}
\end{equation*}
The forward, backward and central difference operators may be defined simply in terms of these shifts as
\begin{equation*}
\delta_{t+} = \frac{1}{h_{t}}\left(e_{t+}-1\right)\qquad\delta_{t-} = \frac{1}{h_{t}}\left(1-e_{t-}\right)\qquad\delta_{to} = \frac{1}{2h_{t}}\left(e_{t+}-e_{t-}\right)
\end{equation*}
(Here, the symbol 1 corresponds to the identity operation.) All of these serve as approximations to a first time derivative; a centered approximation to a second derivative is given by
\begin{equation*}
\delta_{t+}\delta_{t-} = \frac{1}{h_{t}^2}\left(e_{t+}-2+e_{t-}\right)
\end{equation*}

Forward, backward and central time-averaging operators, defined by
\begin{equation*}
\mu_{t+} = \frac{1}{2}\left(e_{t+}+1\right)\qquad\mu_{t-} = \frac{1}{2}\left(1+e_{t-}\right)\qquad\mu_{to} = \frac{1}{2}\left(e_{t+}+e_{t-}\right)
\end{equation*}
are approximations to the identity operator. Another averaging operator, useful in the context of the construction of conservative schemes, is given by
\begin{equation*}
\mu_{t+}\mu_{t-} = \frac{1}{4}\left(e_{t+}+2+e_{t-}\right)
\end{equation*}
Note also that
\begin{eqnarray}
\label{shitid1}
\mu_{t+}\delta_{t-} &=& \mu_{t-}\delta_{t+} = \delta_{to}\\
\label{shitid2}
\mu_{to} &=& 1+\frac{h_{t}^2}{2}\delta_{t+}\delta_{t-}\\
\label{mudid}
\mu_{t-}+\frac{h_{t}}{2}\delta_{t-} &=& 1
\end{eqnarray}

In the energetic analysis of difference schemes, the following identities are indispensable: for any grid function ${\bf f}$, 
\begin{equation}
\label{mudelid}
\left(\mu_{t\star}{\bf f}^{T}\right)\left(\delta_{t\star}{\bf f}\right)= \frac{1}{2}\delta_{t\star}{\bf f}^{T}{\bf f}
\end{equation}
where ``$\star$" stands for any of ``+", ``-" or ``$\cdot$", and
\begin{eqnarray}
\label{shiftid}
{\bf f}^{T}e_{t-}{\bf f} &=& (\mu_{t-}{\bf f}^{T})(\mu_{t-}{\bf f})-\frac{h_{t}^2}{4}(\delta_{t-}{\bf f}^{T})(\delta_{t-}{\bf f})\\
\label{shiftid2}
\mu_{t+}\left({\bf f}^{T}e_{t-}{\bf f}\right) &=& {\bf f}^{T}\mu_{to}{\bf f} 
\end{eqnarray}

The following identities are more useful when examining the conservation of angular momentum. For any two grid functions ${\bf f}$ and ${\bf g}$ of the same number of elements, 
\begin{eqnarray*}
{\bf f}^{T}\delta_{t+}{\bf g} &=& \delta_{t+}\left((\mu_{t-}{\bf f}^{T}){\bf g})\right)-\mu_{t+}\left((\delta_{t-}{\bf f}^{T}){\bf g}\right)\\
(\mu_{to}{\bf f})^{T}\delta_{t+}{\bf g} &=& \delta_{t+}\left((\mu_{t-}{\bf f}^{T}){\bf g})\right)-\frac{1}{2}\left((\delta_{t-}{\bf f}^{T})(e_{t+}{\bf g})+(\delta_{t+}{\bf f}^{T}){\bf g}\right)
\end{eqnarray*}
Both are analogous to the product rule of differentiation. In particular, if for some two-vector grid function ${\bf q}$, it is true that ${\bf f} = \tilde{{\bf q}}$, and ${\bf g} = \delta_{t-}{\bf q}$, it then follows immediately that
\begin{equation}
\label{amid}
\tilde{{\bf q}}^{T}\delta_{t+}\delta_{t-}{\bf q} = (\mu_{to}\tilde{{\bf q}}^{T})\delta_{t+}\delta_{t-}{\bf q} = \delta_{t+}\left((\mu_{t-}\tilde{{\bf q}}^{T})\delta_{t-}{\bf q})\right)
\end{equation}

Another useful identity is the following:
\begin{equation}
\label{tildeid2}
(\mu_{to}\tilde{{\bf f}}^{T}){\bf f} = \frac{h_{t}}{2}\delta_{t+}(\tilde{{\bf f}}^{T}e_{t-}{\bf f})
\end{equation}

\subsection{Spatial Operators}

The only spatial difference operators which will be employed here are the forward and backward difference, defined here in terms of their action on the grid function ${\bf f}_{i}^{n}$ and by
\begin{equation*}
\delta_{x+}{\bf f}_{i}^{n} = \frac{1}{h_{x}}\left({\bf f}_{i+1}^{n}-{\bf f}_{i}^{n}\right)\qquad\delta_{x-}{\bf f}_{i}^{n} = \frac{1}{h_{x}}\left({\bf f}_{i}^{n}-{\bf f}_{i-1}^{n}\right)
\end{equation*}
and are both approximations to a first spatial derivative; a centered approximation to the second derivative is then given by $\delta_{x+}\delta_{x-}$. 

The operators $\delta_{t+}$,  $\delta_{t-}$, $\delta_{to}$, $\mu_{t+}$, $\mu_{t-}$, $\mu_{to}$, $\delta_{x+}$, and $\delta_{x-}$ all commute. 

\subsection{Inner Products and Norms}

The spatial inner product which will be of use in the present article is defined in terms of two grid functions ${\bf f}$ and ${\bf g}$, again of the same number of elements, and over the finite range of indices $i=r,\hdots, s$, for $r$, $s$ integer such that $s\geq r$:
\begin{equation*}
\langle {\bf f}, {\bf g}\rangle_{[r,s]} = h_{x}\sum_{i=r}^{i=s}({\bf f}_{i})^{T}{\bf g}_{i}
\end{equation*}
It is worth noting that for any vector grid functions ${\bf f}$ and ${\bf g}$, and any scalar grid function $l$, 
\begin{equation}
\label{tripleid}
\langle {\bf f}, l{\bf g}\rangle_{[r,s]} = \langle {\bf f}^{T}{\bf g}, l\rangle_{[r,s]}
\end{equation}

The definition of the norm follows in the usual way as
\begin{equation*}
\|{\bf f}\|_{[r,s]} = \langle {\bf f}, {\bf f}\rangle_{[r,s]}^{1/2}
\end{equation*}

The two standard inequalities which follow from the above definitions are the triangle inequality
\begin{equation}
\label{trieq}
\|{\bf f}+{\bf g}\|_{[r,s]} \leq \|{\bf f}\|_{[r,s]}+\|{\bf g}\|_{[r,s]}
\end{equation}
and the Cauchy-Schwartz inequality
\begin{equation*}
|\langle {\bf f}, {\bf g}\rangle_{[r,s]}|\leq \|{\bf f}\|_{[r,s]}\|{\bf g}\|_{[r,s]}
\end{equation*}

Summation by parts follows from the definition of the inner product as
\begin{equation}
\label{ipd}
\langle {\bf f}, \delta_{x+}{\bf g}\rangle_{[r,s]}=-\langle \delta_{x-}{\bf f}, {\bf g}\rangle_{[r+1,s]}+({\bf f}_{s})^{T}{\bf g}_{s+1}-({\bf f}_{r})^{T}{\bf g}_{r}
\end{equation}

It is important to note that the identities defined in Section \ref{timesec} involving the product of two grid functions extend immediately when an inner product is formed. For instance, identity \eqref{tildeid2} becomes
\begin{equation*}
\langle\mu_{to}\tilde{{\bf f}},{\bf f}\rangle_{{\mathcal D}} = \frac{h_{t}}{2}\delta_{t+}\langle \tilde{{\bf f}}, e_{t-}{\bf f}\rangle_{{\mathcal D}}
\end{equation*}
when an inner product is taken over some spatial domain ${\mathcal D}$. 
\subsection{Bounds}

It is straightforward to relate a bound on the norm of a time difference of a grid function to the norm of the grid function itself. For instance, consider a vector-valued grid function ${\bf f}_{i}^{n}$, and suppose it is true that
\begin{equation}
\label{timebound}
\|\delta_{t-}{\bf f}^{n}\|_{{\mathcal D}}\leq K
\end{equation}
over some spatial interval ${\mathcal D}$, for all $n$, and for some constant $K$. It then follows, from the triangle inequality \ref{trieq}, that
\begin{equation}
\label{tdbound}
{\bf f}^{n} = {\bf f}^{0}+h_{t}\sum_{l=1}^{n}\delta_{t-}{\bf f}^{l}\Longrightarrow\|{\bf f}^{n}\|_{{\mathcal D}}\leq \|{\bf f}^{0}\|_{{\mathcal D}}+h_{t}\sum_{l=1}^{n}\|\delta_{t-}{\bf f}^{l}\|_{{\mathcal D}}\leq \|{\bf f}^{0}\|_{{\mathcal D}}+h_{t}nK
\end{equation}
Thus, given a bound such as \eqref{timebound}, growth of the norm of a grid function is at most linear; this is independent of any boundary condition considerations. 

The following bound on the norm of a grid function in terms of its spatial difference follows directly from the triangle inequality \eqref{trieq}:
\begin{equation}
\label{diffbound}
\|\delta_{x-}{\bf f}\|_{[r+1,s]}\leq \frac{2}{h_{x}}\|{\bf f}\|_{[r,s]}
\end{equation}

If, for any scalar grid function $f_{i}$ (which could be a component of a vector grid function), it is true that $f_{r}=0$, then for any $r+1\leq m\leq s$, it is true that
\begin{equation*}
f_{m} = \sum_{i=r+1}^{m}\delta_{x-}f_{i}= \langle 1,\delta_{x-}f\rangle_{[r+1,m]}
\end{equation*}
where ``1" refers to a scalar grid function consisting of a sequence of ones. Then, by the Cauchy-Schwartz inequality, 
\begin{equation*}
|f_{m}|\leq \|1\|_{[r+1,m]}\|\delta_{x-}f\|_{[r+1,m]}\leq \|1\|_{[r+1,s]}\|\delta_{x-}f\|_{[r+1,s]}= \sqrt{h_{x}(s-r)}\|\delta_{x-}f\|_{[r+1,s]}
\end{equation*}
which implies, furthermore, that
\begin{equation}
\label{spatubound}
\|f\|_{[r+1,s]}\leq h_{x}(s-r)\|\delta_{x-}f\|_{[r+1,s]}
\end{equation}

\section{Finite Difference Schemes for Systems {\bf S} and {\bf K}}
\label{sksec}
In Table \ref{tab1}, several finite difference schemes for systems {\bf S} and {\bf K} are presented; these are indicated by {\bf s} and {\bf k}, with a distinguishing superscript. The variety of schemes, in particular for system {\bf S}, is intended to illustrate the many subtle differences among the schemes, with respect to conservation properties, numerical stability, as well as ease of implementation. The number of possible schemes, even of limited stencil, is of course much larger than that indicated here. Typically, in sound synthesis applications, explicit schemes similar to ${\bf s}^{(a)}$ are used \cite{Bankdafx04}.  

\begin{table}
\caption{\label{tab1} Finite difference schemes for systems {\bf S}, and {\bf K}. All instances of a grid function $\xi$ or \etab\, refer to that function at grid location $i$ and time step $n$, i.e., $\xi_{i}^{n}$ and \etab$_{i}^{n}$. The quantities $p$ and ${\bf q}$ are defined in Eq. \eqref{pqdef}. All schemes are defined over the spatial interval $i\in{\mathcal D}$, and for $n\geq 0$. The set ${\mathcal D}^{+}$ is defined in the third paragraph of Section \ref{sksec}.} 
\begin{center}\begin{scriptsize}
\begin{tabular}{l|c|l}
&\mbox{Defining Equations} \\\hline\hline
\raisebox{-0.1in}{${\bf s}^{(a)}$} && $\phi_{{\bf s}^{(a)}} = p + \frac{1-\alpha}{2}{\bf q}^{T}{\bf q}$  \\
& & $\psib_{{\bf s}^{(a)}} = \alpha{\bf q}+ \frac{1-\alpha}{2}({\bf q}^{T}{\bf q}+2p){\bf q}$  \\\cline{3-3}\cline{1-1} 
\raisebox{-0.1in}{${\bf s}^{(b)}$} && $\phi_{{\bf s}^{(b)}} = p + \frac{1-\alpha}{2}{\bf q}^{T}{\bf q}$   \\
&& $\psib_{{\bf s}^{(b)}} = \alpha{\bf q}+ \frac{1-\alpha}{2}({\bf q}^T{\bf q}+2p)(\mu_{to}{\bf q})$\\\cline{3-3}\cline{1-1} 
\raisebox{-0.1in}{${\bf s}^{(c)}$} & $\delta_{t+}\delta_{t-}\xi = \delta_{x+}\phi_{\bullet}$& $\phi_{{\bf s}^{(c)}}=p + \frac{1-\alpha}{2}{\bf q}^{T}\mu_{to}{\bf q}$   \\
& $\delta_{t+}\delta_{t-}\etab = \delta_{x+}\psib_{\bullet}$ & $\psib_{{\bf s}^{(c)}} = \alpha{\bf q}+ \frac{1-\alpha}{2}{\bf q}^T{\bf q}\mu_{to}{\bf q}+2(\mu_{t+}\mu_{t-}p)({\bf q})$\\\cline{3-3}\cline{1-1} \raisebox{-0.1in}{${\bf s}^{(d)}$} && $ \phi_{{\bf s}^{(d)}}= p + \frac{1-\alpha}{2}{\bf q}^{T}\mu_{to}{\bf q}$   \\
&& $ \psib_{{\bf s}^{(d)}}= \alpha{\bf q}+ \frac{1-\alpha}{2}{\bf q}^T(\mu_{to}{\bf q}){\bf q}+2(\mu_{t+}\mu_{t-}p)({\bf q})$\\\cline{3-3}\cline{1-1} 
\raisebox{-0.1in}{${\bf s}^{(e)}$} && $ \phi_{{\bf s}^{(e)}}= p + \frac{1-\alpha}{2}\mu_{to}({\bf q}^{T}{\bf q})$   \\
&& $ \psib_{{\bf s}^{(e)}}= \alpha{\bf q}+ \frac{1-\alpha}{2}(\mu_{to}({\bf q}^T{\bf q}+2p))(\mu_{to}{\bf q})$\\\hline ${\bf k}^{(a)}$ & \raisebox{-0.1in}{$\delta_{t+}\delta_{t-}\etab = \alpha{\mathcal G}_{\bullet}\delta_{x+}{\bf q}$} & ${\mathcal G}_{{\bf k}^{(a)}} = 1+\frac{1}{2\alpha}\|{\bf q}\|_{D^{+}}^{2}$\\\cline{3-3}\cline{1-1}
${\bf k}^{(b)}$ && ${\mathcal G}_{{\bf k}^{(b)}} = 1+\frac{1}{2\alpha}\mu_{t+}\langle{\bf q}, e_{t-}{\bf q}\rangle_{D^{+}}$\\
 \end{tabular}
\end{scriptsize}\end{center}
\end{table}

All the schemes given in the table for system {\bf S} have the same form, indicated in the second column of the table, which is a direct discretization of Eqs. \eqref{Morse3}, in $\phi$ and $\psib$. The shorthand forms
\begin{equation}
\label{pqdef}
p = \delta_{x-}\xi\qquad {\bf q} = \delta_{x-}\etab
\end{equation}
are used throughout the rest of this article. Distinctions among the various schemes are due to variations in the discretization of $\phi$ and $\psib$, given explicitly in the third column of the table. Similarly, for system {\bf K}, the two schemes given have the form of a direct discretization of Eq. \eqref{kc}, and variations are due to the way in which the quantity ${\mathcal G}$ is discretized. 

All the schemes below are consistent with systems {\bf S} or {\bf K} and accurate to second order in both time and space (it is simplest to see this by rewriting the schemes in a first-order transmission-line form, as per \cite{Bilbaojasa05}). They are two-step schemes, and for initialization, values of the grid functions $\xi^n$ and $\etab^n$ are required at the first two time steps, i.e., for $n=0$ and $n=1$. The spatial domain of the problem will be limited to $i\in{\mathcal D} = [0,\hdots,N]$. In dealing with boundary conditions, the set ${\mathcal D}^{+} = [1,\hdots,N]$ is also of great utility. 

As mentioned above, each scheme can be considered from various points of view. Does it possess conserved analogues of angular momentum and energy? Are there simple conditions under which numerical stability can be ensured? Is the scheme explicit or implicit, and if implicit, do existence and uniqueness conditions for the numerical solution follow? All of these points will be dealt with in turn in the following sections.

\subsection{Conservation of Angular Momentum}

The conservation of angular momentum (in the $x$-direction) is perhaps the simplest property to examine. Consider first the simple scheme {\bf s}$^{(a)}$, as given in Table \ref{tab1}. Taking the inner product of the second equation with $\tilde{\etab}$ over the domain ${\mathcal D}$ gives 
\begin{eqnarray*}
\langle \tilde{\etab}, \delta_{t+}\delta_{t-}\etab\rangle_{{\mathcal D}} &=& \langle \tilde{\etab},\delta_{x-}\psib_{{\bf s}^{(a)}}\rangle_{{\mathcal D}}\\
&=& -\langle \tilde{{\bf q}},\psib_{{\bf s}^{(a)}}\rangle_{{\mathcal D}^{+}} + \mathcal{B}_{{\mathcal A}, {\bf s}^{(a)}}\\
&=& \mathcal{B}_{{\mathcal A}, {\bf s}^{(a)}}
\end{eqnarray*}
where the second and third equalities above follow from summation by parts (Eq. \eqref{ipd}) and identity \eqref{tildeid} (note that $\psib_{{\bf s}^{(a)}}$ is a scalar multiple of ${\bf q}$). The boundary term $\mathcal{B}_{{\mathcal A}, {\bf s}^{(a)}}$ is given in Table \ref{tab2}. From identity \eqref{amid}, it then follows that
\begin{equation*}
\delta_{t+}{\mathcal A}_{{\bf s}^{(a)}} = \mathcal{B}_{{\mathcal A}, {\bf s}^{(a)}}
\end{equation*}
where 
\begin{equation*}
{\mathcal A}_{{\bf s}^{(a)}} = \langle \mu_{t-}\tilde{\etab}, \delta_{t-}\etab\rangle_{{\mathcal D}}
\end{equation*}
can be identified easily with the angular momentum. Clearly, if the boundary term $\mathcal{B}_{{\mathcal A}, {\bf s}^{(a)}}$ vanishes, then angular momentum is conserved by scheme ${\bf s}^{(a)}$.

 \begin{table}
\caption{\label{tab2} Conserved angular momentum, and boundary terms for the schemes given in Table \ref{tab1}. In all cases for which expressions are provided, it is true that $\delta_{t+}{\mathcal A}_{\bullet} = {\mathcal B}_{\bullet}$. The lack of a conserved angular momentum is indicated by the symbol ``$- -$". } 
\begin{center}\begin{scriptsize}
\begin{tabular}{l|l|l}
&\mbox{Conserved angular momentum} & Boundary term\\\hline\hline
${\bf s}^{(a)}$ & ${\mathcal A}_{{\bf s}^{(a)}} = \langle \mu_{t-}\tilde{\etab}, \delta_{t-}\etab\rangle_{{\mathcal D}}$& ${\mathcal B}_{{\mathcal A}, {\bf s}^{(a)}} = \tilde{\etab}^{T}_{N}\psib_{N+1}-\tilde{\etab}^{T}_{0}\psib_{0}$\\\hline
${\bf s}^{(b)}$ & ${\mathcal A}_{{\bf s}^{(b)}} = \langle \mu_{t-}\tilde{\etab}, \delta_{t-}\etab\rangle_{{\mathcal D}}+ \frac{\alpha h_{t}}{2}\langle \tilde{{\bf q}}, e_{t-}{\bf q}\rangle_{{\mathcal D}^{+}}$& ${\mathcal B}_{{\mathcal A}, {\bf s}^{(b)}} = (\mu_{to}\tilde{\etab}^{T}_{N})\psib_{N+1}-(\mu_{to}\tilde{\etab}^{T}_{0})\psib_{0}$  \\\hline
${\bf s}^{(c)}$ & $- -$& $- -$  \\\hline
${\bf s}^{(d)}$ & ${\mathcal A}_{{\bf s}^{(d)}} = \langle \mu_{t-}\tilde{\etab}, \delta_{t-}\etab\rangle_{{\mathcal D}}$& ${\mathcal B}_{{\mathcal A}, {\bf s}^{(d)}} = \tilde{\etab}^{T}_{N}\psib_{N+1}-\tilde{\etab}^{T}_{0}\psib_{0}$  \\\hline
${\bf s}^{(e)}$ & ${\mathcal A}_{{\bf s}^{(e)}} = \langle \mu_{t-}\tilde{\etab}, \delta_{t-}\etab\rangle_{{\mathcal D}}+ \frac{\alpha h_{t}}{2}\langle \tilde{{\bf q}}, e_{t-}{\bf q}\rangle_{{\mathcal D}^{+}}$& ${\mathcal B}_{{\mathcal A}, {\bf s}^{(e)}} = (\mu_{to}\tilde{\etab}^{T}_{N})\psib_{N+1}-(\mu_{to}\tilde{\etab}^{T}_{0})\psib_{0}$ \\\hline
${\bf k}^{(a)}$ & ${\mathcal A}_{{\bf k}^{(a)}} = \langle \mu_{t-}\tilde{\etab}, \delta_{t-}\etab\rangle_{{\mathcal D}}$& ${\mathcal B}_{{\mathcal A}, {\bf k}^{(a)}} = {\mathcal G}_{{\bf k}^{(a)}}\left(\tilde{\etab}^{T}_{N}{\bf q}_{N+1}-\tilde{\etab}^{T}_{0}{\bf q}_{0}\right)$  \\\hline
${\bf k}^{(b)}$ & ${\mathcal A}_{{\bf k}^{(b)}} = \langle \mu_{t-}\tilde{\etab}, \delta_{t-}\etab\rangle_{{\mathcal D}}$& ${\mathcal B}_{{\mathcal A}, {\bf k}^{(b)}} = {\mathcal G}_{{\bf k}^{(b)}}\left(\tilde{\etab}^{T}_{N}{\bf q}_{N+1}-\tilde{\etab}^{T}_{0}{\bf q}_{0}\right)$  \\\hline
\end{tabular}
\end{scriptsize}\end{center}
\end{table}

The treatment of system ${\bf s}^{(b)}$ is similar, except that it is now necessary to take the inner product of the second equation with $\mu_{to}\tilde{\etab}$, instead of $\tilde{\etab}$, giving
\begin{equation*}
 \langle \mu_{to}\tilde{\etab}, \delta_{t+}\delta_{t-}\etab\rangle_{{\mathcal D}} = \langle \mu_{to}\tilde{\etab},\delta_{x-}\psib_{{\bf s}^{(b)}}\rangle_{{\mathcal D}}= -\langle \mu_{to}\tilde{{\bf q}},\psib_{{\bf s}^{(b)}}\rangle_{{\mathcal D}^{+}} + \mathcal{B}_{{\mathcal A}, {\bf s}^{(b)}}
\end{equation*}
where the boundary term $\mathcal{B}_{{\mathcal A}, {\bf s}^{(b)}}$ is given in Table \ref{tab2}. The expression for $\psib_{{\bf s}^{(b)}}$, from Table \ref{tab1}, is made up of two terms; the first is a scalar multiple of ${\bf q}$, and the second a multiple of $\mu_{to}{\bf q}$. Thus, applying identity \eqref{tildeid}, 
\begin{eqnarray*}
\langle \mu_{to}\tilde{\etab}, \delta_{t+}\delta_{t-}\etab\rangle_{{\mathcal D}} &=& -\alpha\langle \delta_{x-}\mu_{to}\tilde{\etab},{\bf q}\rangle_{{\mathcal D}^{+}} + \mathcal{B}_{{\mathcal A}, {\bf s}^{(b)}}\\
&=& -\alpha\langle \mu_{to}\tilde{{\bf q}},{\bf q}\rangle_{{\mathcal D}^{+}} + \mathcal{B}_{{\mathcal A}, {\bf s}^{(b)}}\\
&=& \frac{-\alpha h_{t}}{2}\delta_{t+}\langle\tilde{{\bf q}},e_{t-}{\bf q}\rangle_{{\mathcal D}^{+}} + \mathcal{B}_{{\mathcal A}, {\bf s}^{(b)}}
\end{eqnarray*}
where, in the second and third equalities above, commutativity of the operators $\mu_{to}$ and $\delta_{x-}$ and identity \eqref{tildeid2} have been used, respectively. 
Finally, applying identity \eqref{tildeid}, this can be rewritten as
\begin{equation*}
\delta_{t+}{\mathcal A}_{{\bf s}^{(b)}} = \mathcal{B}_{{\mathcal A}, {\bf s}^{(b)}}
\end{equation*}
where 
\begin{equation*}
{\mathcal A}_{{\bf s}^{(b)}} = \langle \mu_{t-}\tilde{\etab}, \delta_{t-}\etab\rangle_{{\mathcal D}}+ \frac{\alpha h_{t}}{2}\langle\tilde{{\bf q}}, e_{t-}{\bf q}\rangle_{{\mathcal D}^{+}}
\end{equation*}

This form of the angular momentum is distinct from the quantity conserved under scheme ${\bf s}^{(a)}$, but note that in the limit as $h_{t}$ becomes small, the two definitions approach one another. 
Scheme ${\bf s}^{(c)}$ does not possess a simple conserved quantity analogous to an angular momentum, but schemes ${\bf s}^{(d)}$ and ${\bf s}^{(e)}$ do---their analysis is nearly identical to that of ${\bf s}^{(a)}$ and ${\bf s}^{(b)}$, respectively, and the conserved quantities and boundary terms are given in Table \ref{tab2}. 

For schemes ${\bf k}^{(a)}$ and ${\bf k}^{(b)}$, the analysis is very similar, and conservation of angular momentum holds in either case, with conserved quantities as given in Table \ref{tab2}.

\subsection{Conservation of Energy}
\label{egyconssec}
Conservation of energy for the schemes ${\bf s}^{(\bullet)}$ is always arrived at in the following way: given the general form of the scheme, shown in the second column of Table \ref{tab1}, take the inner product over the domain ${\mathcal D}$ of the first equation with $\delta_{to}\xi$, and the second with $\delta_{to}\etab$. After using summation by parts (Eq. \eqref{ipd}) and adding the resulting equations one arrives at
\begin{equation}
\label{egybal1}
\delta_{t+}\Big[\frac{1}{2}\left(\|\delta_{t-}\xi\|_{{\mathcal D}}^{2}+\|\delta_{t-}\xi\|_{{\mathcal D}}^{2}\right)\Big] +\langle \delta_{to} p,\phi_{{\bf s}^{(\bullet)}}\rangle_{{\mathcal D}^{+}}+\langle \delta_{to} {\bf q},\psib_{{\bf s}^{(\bullet)}}\rangle_{{\mathcal D}^{+}}={\mathcal B}_{{\mathcal H}, {\bf s}^{(\bullet)}}
\end{equation}
where $\phi_{{\bf s}^{(\bullet)}}$ and $\psib_{{\bf s}^{(\bullet)}}$ depend on the choice of scheme, and are given in Table \ref{tab1}. The boundary term ${\mathcal B}_{{\mathcal H},{\bf s}^{(\bullet)}}$ is of the same form for all the schemes given for system {\bf S}:
\begin{equation}
\label{bhdef}
{\mathcal B}_{{\mathcal H},{\bf s}^{(\bullet)}} = (\delta_{to}\xi_{N})\phi_{{\bf s}^{(\bullet)},N+1}+(\delta_{to}\etab_{N}^{T})\psib_{{\bf s}^{(\bullet)},N+1}-(\delta_{to}\xi_{0})\phi_{{\bf s}^{(\bullet)},0}-(\delta_{to}\etab_{0}^{T})\psib_{{\bf s}^{(\bullet)},0}
\end{equation}
Clearly, the first term on the left hand side of Eq. \eqref{egybal1} behaves as a difference approximation to the first derivative of the kinetic energy. The second two terms will vary from scheme to scheme, depending on the forms of $\phi_{{\bf s}^{(\bullet)}}$ and $\psib_{{\bf s}^{(\bullet)}}$. 

Scheme ${\bf s}^{(a)}$ does not possess a conserved energy, but the other four do. For example, for scheme ${\bf s}^{(b)}$, using the forms of $\phi_{{\bf s}^{(b)}}$ and $\psib_{{\bf s}^{(b)}}$ given in Table \ref{tab1}, one may write
\begin{eqnarray}
\label{crap1}
\langle\delta_{to} p,\phi_{{\bf s}^{(b)}}\rangle_{{\mathcal D}^{+}} &=&\langle\delta_{to} p,p + \frac{1-\alpha}{2}{\bf q}^{T}{\bf q}\rangle_{{\mathcal D}^{+}}\\
&=&\langle\delta_{to}p,\alpha p + \frac{1-\alpha}{2}\left({\bf q}^{T}{\bf q}+2p\right)\rangle_{{\mathcal D}^{+}}\notag\\
&=& \delta_{t+}\Big[\frac{\alpha}{2}\langle p, e_{t-}p\rangle_{{\mathcal D}^{+}}\Big]+\frac{1-\alpha}{2}\langle \delta_{to} p,{\bf q}^{T}{\bf q}+2p\rangle_{{\mathcal D}^{+}}\notag 
\end{eqnarray}
and
\begin{eqnarray}
\label{crap2}
\langle\delta_{to} {\bf q},\psib_{{\bf s}^{(b)}}\rangle_{{\mathcal D}^{+}} &=&\langle\delta_{to}{\bf q},\alpha{\bf q}+ \frac{1-\alpha}{2}({\bf q}^T{\bf q}+2p)(\mu_{to}{\bf q})\rangle_{{\mathcal D}^{+}}\\
&=& \delta_{t+}\Big[\frac{\alpha}{2}\langle q, e_{t-}q\rangle_{{\mathcal D}^{+}}\Big]+\frac{1-\alpha}{2}\langle \delta_{to}{\bf q},({\bf q}^T{\bf q}+2p)(\mu_{to}{\bf q})\rangle_{{\mathcal D}^{+}}\notag\\ 
&=& \delta_{t+}\Big[\frac{\alpha}{2}\langle q, e_{t-}q\rangle_{{\mathcal D}^{+}}\Big]+\frac{1-\alpha}{4}\langle \delta_{to}({\bf q}^{T}{\bf q}),({\bf q}^T{\bf q}+2p)\rangle_{{\mathcal D}^{+}}\notag
\end{eqnarray}
and thus, combining Eqs. \eqref{crap1} and \eqref{crap2},
\begin{eqnarray*}
\langle \delta_{to} p,\phi_{{\bf s}^{(b)}}\rangle_{{\mathcal D}^{+}}+ \langle \delta_{to} {\bf q},\psib_{{\bf s}^{(b)}}\rangle_{{\mathcal D}^{+}} &=& \delta_{t+}\Big[\frac{\alpha}{2}\langle p, e_{t-}p\rangle_{{\mathcal D}^{+}} + \frac{\alpha}{2}\langle q, e_{t-}q\rangle_{{\mathcal D}^{+}}\Big]\\
&&+ \frac{1-\alpha}{2}\langle \delta_{to}(\frac{1}{2}{\bf q}^{T}{\bf q}+p),\frac{1}{2}{\bf q}^T{\bf q}+2p\rangle_{{\mathcal D}^{+}}\notag\\
&=& \delta_{t+}{\mathcal V}_{{\bf s}^{(b)}}\notag
\end{eqnarray*}
where ${\mathcal V}_{{\bf s}^{(b)}}$ is as given in Table \ref{tab3}. Thus, one has, immediately, 
\begin{equation*}
\delta_{t+}{\mathcal H}_{{\bf s}^{(b)}} = \delta_{t+}\left({\mathcal T}_{{\bf s}^{(b)}}+{\mathcal V}_{{\bf s}^{(b)}}\right) = {\mathcal B}_{{\mathcal H}, {\bf s}^{(b)}}
\end{equation*}

\begin{table}
\caption{\label{tab3} Discrete kinetic and potential energies for the schemes given in Table \ref{tab1}; their sum will be conserved. The symbol $- -$ indicates that the scheme is not conservative. When such quantities are indicated, the discrete energy conservation property $\delta_{t+}\left({\mathcal T}_{{\bf s}^{(\bullet)}}+{\mathcal V}_{{\bf s}^{(\bullet)}}\right) = {\mathcal B}_{{\mathcal H}, {\bf s}^{(\bullet)}}$ holds, where ${\mathcal B}_{{\mathcal H}, {\bf s}^{(\bullet)}}$ is given in Eq. \eqref{bhdef}.  }
\begin{center}\begin{scriptsize}
\begin{tabular}{l|l|l}
& Kinetic energy & Potential energy\\\hline\hline
${\bf s}^{(a)}$ & $- -$ & $- -$  \\\hline
\raisebox{-0.1in}{${\bf s}^{(b)}$} & & ${\mathcal V}_{{\bf s}^{(b)}} = \frac{\alpha}{2}\langle p, e_{t-}p\rangle_{{\mathcal D}^{+}}+\frac{\alpha}{2}\langle {\bf q}, e_{t-}{\bf q}\rangle_{{\mathcal D}^{+}}$  \\
&&\qquad\qquad$+\frac{1-\alpha}{2}\langle p+\frac{1}{2}{\bf q}^{T}{\bf q}, e_{t-}\left(p+\frac{1}{2}{\bf q}^{T}{\bf q}\right)\rangle_{{\mathcal D}^{+}}$\\\cline{3-3}\cline{1-1}
\raisebox{-0.1in}{${\bf s}^{(c)}$} & & ${\mathcal V}_{{\bf s}^{(c)}} = \frac{1}{2}\langle p, e_{t-}p\rangle_{{\mathcal D}^{+}}+\frac{\alpha}{2}\langle {\bf q}, e_{t-}{\bf q}\rangle_{{\mathcal D}^{+}}$  \\
&&\qquad\qquad$+\frac{1-\alpha}{2}\left(\|\mu_{t-}p+\frac{1}{2}{\bf q}^{T}e_{t-}{\bf q}\|^{2}_{{\mathcal D}^{+}}-\|\mu_{t-}p\|_{{\mathcal D}^{+}}^2\right)$\\
&${\mathcal T}_{{\bf s}^{(\bullet)}} = \frac{1}{2}\left(\|\delta_{t-}\xi\|_{{\mathcal D}}^{2}+\|\delta_{t-}\etab\|_{{\mathcal D}}^{2}\right)$& \quad\quad\quad$\quad+\frac{1-\alpha}{8}\|\tilde{{\bf q}}^{T}e_{t-}{\bf q}\|_{{\mathcal D}^{+}}^2$\\\cline{3-3}\cline{1-1}
\raisebox{-0.1in}{${\bf s}^{(d)}$} & & ${\mathcal V}_{{\bf s}^{(d)}} = \frac{1}{2}\langle p, e_{t-}p\rangle_{{\mathcal D}^{+}}+\frac{\alpha}{2}\langle {\bf q}, e_{t-}{\bf q}\rangle_{{\mathcal D}^{+}}$ \\
&& \qquad\qquad$+\frac{1-\alpha}{2}\left(\|\mu_{t-}p+\frac{1}{2}{\bf q}^{T}e_{t-}{\bf q}\|^{2}_{{\mathcal D}^{+}}-\|\mu_{t-}p\|_{{\mathcal D}^{+}}^2\right)$ \\\cline{3-3}\cline{1-1}
\raisebox{-0.1in}{${\bf s}^{(e)}$} & & ${\mathcal V}_{{\bf s}^{(e)}} = \frac{1}{2}\langle p, e_{t-}p\rangle_{{\mathcal D}^{+}}+\frac{\alpha}{2}\langle {\bf q}, e_{t-}{\bf q}\rangle_{{\mathcal D}^{+}}$ \\
&&\qquad\qquad$+\frac{1-\alpha}{2}\mu_{t-}\left(\|p+\frac{1}{2}{\bf q}^{T}{\bf q}\|^{2}_{{\mathcal D}^{+}}-\|p\|_{{\mathcal D}^{+}}^2\right)$\\\hline
${\bf k}^{(a)}$ & $- -$ & $- -$ \\\hline
${\bf k}^{(b)}$ & ${\mathcal T}_{{\bf k}^{(b)}} = \frac{1}{2}\|\delta_{t-}\etab\|_{{\mathcal D}}^{2}$& ${\mathcal V}_{{\bf k}^{(b)}} = \frac{\alpha}{2}\langle {\bf q}, e_{t-}{\bf q}\rangle_{{\mathcal D}^{+}}\left(1+\frac{1}{4\alpha}\langle {\bf q}, e_{t-}{\bf q}\rangle_{{\mathcal D}^{+}}\right)$ \\
\end{tabular}
\end{scriptsize}\end{center}
\end{table}

One may proceed in a similar vein for schemes ${\bf s}^{(c)}$, ${\bf s}^{(d)}$ and ${\bf s}^{(e)}$; all these algorithms are energy-conserving, with expressions for discrete kinetic and potential energy given in Table \ref{tab3}. 

The two schemes for system {\bf K}, ${\bf k}^{(a)}$ and ${\bf k}^{(b)}$, given in Table \ref{tab1} differ only in the treatment of the quantity ${\mathcal G}$. For either scheme, one may take the inner product with $\delta_{to}\etab$ to get
\begin{equation*}
\frac{1}{2}\langle\delta_{to}\etab, \delta_{t+}\delta_{t-}\etab\rangle_{{\mathcal D}} = \alpha{\mathcal G}_{{\bf k}^{(\bullet)}}\langle\delta_{to}\etab, \delta_{x+}{\bf q}\rangle_{{\mathcal D}}
\end{equation*}
and, again using summation by parts (Eq. \eqref{ipd}), 
\begin{eqnarray*}
\delta_{t+}\Big[\frac{1}{2}\|\delta_{t-}\etab\|_{{\mathcal D}}^{2}\Big]+\frac{\alpha}{2}{\mathcal G}_{{\bf k}^{(\bullet)}}\delta_{t+}\langle{\bf q},e_{t-}{\bf q}\rangle_{{\mathcal D}^{+}}&=&{\mathcal B}_{{\mathcal H}, {\bf k}^{(\bullet)}}
\end{eqnarray*}
Here, the boundary term ${\mathcal B}_{{\mathcal H}, {\bf k}^{(\bullet)}}$ is given by
\begin{equation}
\label{bhkdef}
{\mathcal B}_{{\mathcal H}, {\bf k}^{(\bullet)}} = \alpha{\mathcal G}_{{\bf k}^{(\bullet)}}\left((\mu_{t+}\delta_{t-}\etab^{T}_{N}){\bf q}_{N+1}-(\mu_{t+}\delta_{t-}\etab^{T}_{0}){\bf q}_{0}\right)
\end{equation}
For scheme ${\bf k}^{(a)}$, substitution of the expression ${\mathcal G}_{{\bf k}^{(a)}}$ (given in Table \ref{tab1}) does not lead to an energy-conservation property. But for scheme ${\bf k}^{(b)}$, using ${\mathcal G}_{{\bf k}^{(b)}}$ (also given in Table \ref{tab1}) yields
\begin{equation*}
\delta_{t+}\Big[\frac{1}{2}\|\delta_{t-}\etab\|_{{\mathcal D}}^{2}\Big]+\frac{\alpha}{2}\left(1+\frac{1}{2\alpha}\mu_{t+}\langle{\bf q},e_{t-}{\bf q}\rangle_{{\mathcal D}^{+}}\right)\delta_{t+}\langle{\bf q},e_{t-}{\bf q}\rangle_{{\mathcal D}^{+}}={\mathcal B}_{{\mathcal H}, {\bf k}^{(\bullet)}}
\end{equation*}
and, using identity \eqref{mudelid}, 
\begin{equation*}
\delta_{t+}\Big[\frac{1}{2}\|\delta_{t-}\etab\|_{{\mathcal D}}^{2}\Big]+\frac{\alpha}{2}\delta_{t+}\langle{\bf q},e_{t-}{\bf q}\rangle_{{\mathcal D}^{+}}+\frac{1}{8}\delta_{t+}\left(\langle{\bf q},e_{t-}{\bf q}\rangle_{{\mathcal D}^{+}}^{2}\right)={\mathcal B}_{{\mathcal H}, {\bf k}^{(\bullet)}}
\end{equation*}
Finally, one arrives at
\begin{equation*}
\delta_{t+}{\mathcal H}_{{\bf k}^{(b)}} = \delta_{t+}\left({\mathcal T}_{{\bf k}^{(b)}}+{\mathcal V}_{{\bf k}^{(b)}}\right) = {\mathcal B}_{{\mathcal H}, {\bf k}^{(b)}}
\end{equation*}
where the expressions ${\mathcal T}_{{\bf k}^{(b)}}$, and ${\mathcal V}_{{\bf k}^{(b)}}$ are as given in Table \ref{tab3}. This is the desired discrete energy conservation property.

\subsection{Conservative Boundary Conditions}
\label{bcdsec}
Fixed boundary conditions, defined by
\begin{equation}
\label{bcdef1fd}
\xi = 0\qquad \etab = {\bf 0}
\end{equation}
are a direct counterpart to the continuous conditions \eqref{bcdef1}, and are assumed to hold at an endpoint of the interval ${\mathcal D}$, i.e., for $i=0$ or $i=N$. 

Free boundary conditions \eqref{bcdef2} at, e.g., the left end of the discrete domain, can be approximated by
\begin{equation}
\label{bcdef2fda}
p_{0} = 0\qquad {\bf q}_{0} = {\bf 0}
\end{equation}
where it is recalled that for all the schemes discussed here, the definitions \eqref{pqdef} hold. At the right end, such free conditions are given by
\begin{equation}
\label{bcdef2fdb}
p_{N+1} = 0\qquad {\bf q}_{N+1} = {\bf 0}
\end{equation}
For systems ${\bf k}^{(\bullet)}$, which depend only on \etab, the second of each pair of conditions above  suffices to characterize a boundary condition as fixed or free. 

Given the forms of ${\mathcal B}_{{\mathcal A}, {\bf s}^{(\bullet)}}$, ${\mathcal B}_{{\mathcal A}, {\bf k}^{(\bullet)}}$ given in Table \ref{tab2}, and ${\mathcal B}_{{\mathcal H}, {\bf s}^{(\bullet)}}$ and ${\mathcal B}_{{\mathcal H}, {\bf k}^{(\bullet)}}$ given by definitions \eqref{bhdef} and \eqref{bhkdef}, respectively, and recalling the various forms of $\phi_{{\bf s}^{(\bullet)}}$ and $\psib_{{\bf s}^{(\bullet)}}$ given in Table \ref{tab1}, it should be clear that a choice of a fixed discrete boundary condition \eqref{bcdef1fd}, or a free boundary condition such as \eqref{bcdef2fda} or \eqref{bcdef2fdb} at each of the endpoints of the domain ${\mathcal D}$ leads to a vanishing of the boundary terms. Under such conditions, any such scheme is fully conservative, i.e., one has

\begin{equation*}
\delta_{t+}{\mathcal A}_{\bullet} = 0\qquad\delta_{t+}{\mathcal H}_{\bullet} = 0
\end{equation*}
and thus
\begin{equation*}
{\mathcal A}_{\bullet}^n = {\mathcal A}_{\bullet}^0\qquad{\mathcal H}_{\bullet}^n = {\mathcal H}_{\bullet}^0 
\end{equation*}
Such conservative boundary conditions will be assumed for the remainder of this article.

\subsection{Numerical Stability}

Under further conditions, numerical stability of the difference schemes may follow immediately from the discrete conserved energy quantities given in Table \ref{tab3}. (The schemes which do not possess an energy conservation property, namely ${\bf s}^{(a)}$ and ${\bf k}^{(a)}$, will not be examined in this section.) The main goal, in this section, is to find conditions, if they exist, under which the discrete conserved energy is positive for all possible choices of the state variables $\xi$ and $\etab$. If such conditions exist, then bounds on the solution size in terms of initial conditions (and hence a numerical stability guarantee) can be obtained. 

\subsubsection{Positivity of Discrete Conserved Energy}

For the schemes ${\bf s}^{(b)}$, ${\bf s}^{(c)}$, ${\bf s}^{(d)}$, and ${\bf s}^{(e)}$, a discrete energy conservation property exists, and in particular, under conservative boundary conditions such as those discussed in Section \ref{bcdsec}, one has
\begin{equation*}
{\mathcal T}_{{\bf s}^{\bullet}}^{n} + {\mathcal V}_{{\bf s}^{\bullet}}^{n} = {\mathcal H}_{{\bf s}^{\bullet}}^{n} = {\mathcal H}_{{\bf s}^{\bullet}}^{0}
\end{equation*}

It is useful to begin with scheme ${\bf s}^{(d)}$, which has a conserved energy of a particularly simple form. Considering the form for ${\mathcal V}_{{\bf s}^{(d)}}$ given in Table \ref{tab3}, which is 
\begin{eqnarray*}
{\mathcal V}_{{\bf s}^{(d)}} &=& \frac{1}{2}\langle p, e_{t-}p\rangle_{{\mathcal D}^{+}}+\frac{\alpha}{2}\langle {\bf q}, e_{t-}{\bf q}\rangle_{{\mathcal D}^{+}}\\
\qquad &&+\frac{1-\alpha}{2}\left(\|\mu_{t-}p+\frac{1}{2}{\bf q}^{T}e_{t-}{\bf q}\|^{2}_{{\mathcal D}^{+}}-\|\mu_{t-}p\|_{{\mathcal D}^{+}}^2\right)
\end{eqnarray*}
one may then write, employing identity \eqref{shiftid}, 
\begin{eqnarray}
\label{vd}
{\mathcal V}_{{\bf s}^{(d)}} &=& \frac{\alpha}{2}\left(\|\mu_{t-}p\|^2_{{\mathcal D}^{+}}+ \|\mu_{t-}{\bf q}\|^2_{{\mathcal D}^{+}}\right)+\frac{1-\alpha}{2}\left(\|\mu_{t-}p+\frac{1}{2}{\bf q}^{T}e_{t-}{\bf q}\|^{2}_{{\mathcal D}^{+}}\right)\\
\qquad&& -\frac{h_{t}^2}{8}\|\delta_{t-}p\|^2_{{\mathcal D}^{+}}-\frac{\alpha h_{t}^2}{8}\|\delta_{t-}{\bf q}\|^2_{{\mathcal D}^{+}}\\
\end{eqnarray}
Applying definitions \eqref{pqdef} and the bound \eqref{diffbound} yields the following inequality
\begin{eqnarray*}
{\mathcal V}_{{\bf s}^{(d)}}&\geq&\frac{\alpha}{2}\left(\|\mu_{t-}p\|^2_{{\mathcal D}^{+}}+ \|\mu_{t-}{\bf q}\|^2_{{\mathcal D}^{+}}\right)+\frac{1-\alpha}{2}\left(\|\mu_{t-}p+\frac{1}{2}{\bf q}^{T}e_{t-}{\bf q}\|^{2}_{{\mathcal D}^{+}}\right)\\
\qquad&& -\frac{\lambda^2}{2}\|\delta_{t-}\xi\|^2_{{\mathcal D}}-\frac{\alpha\lambda^2}{2}\|\delta_{t-}\etab\|^2_{{\mathcal D}}
\end{eqnarray*}
Note that the parameter $\lambda$, defined in Eq. \eqref{Courantdef}, has been introduced here. One further has
\begin{eqnarray}
\label{sdegy}
{\mathcal H}_{{\bf s}^{(d)}}={\mathcal T}_{{\bf s}^{(d)}}+{\mathcal V}_{{\bf s}^{(d)}}&\geq& \frac{\alpha}{2}\left(\|\mu_{t-}p\|^2_{{\mathcal D}^{+}}+ \|\mu_{t-}{\bf q}\|^2_{{\mathcal D}^{+}}\right)\\
\qquad&&+\frac{1-\alpha}{2}\left(\|\mu_{t-}p+\frac{1}{2}{\bf q}^{T}e_{t-}{\bf q}\|^{2}_{{\mathcal D}^{+}}\right)\notag\\
\qquad&&+ \left(\frac{1}{2}-\frac{\lambda^2}{2}\right)\|\delta_{t-}\xi\|^2_{{\mathcal D}}+\left(\frac{1}{2}-\frac{\alpha\lambda^2}{2}\right)\|\delta_{t-}\etab\|^2_{{\mathcal D}}\notag
\end{eqnarray}

Keeping in mind the condition \eqref{EATcond}, which is assumed a priori, then the discrete conserved energy ${\mathcal H}_{{\bf s}^{(d)}}$ will be non-negative under the conditions
\begin{subequations}
\label{CFLcond}
\begin{eqnarray}
\label{CFLcond1}
\lambda&\leq& \sqrt{1/\alpha}\\
\label{CFLcond2}
\lambda&\leq& 1
\end{eqnarray}
\end{subequations}
which have the form of Courant-Friedrichs-Lewy type conditions \cite{Strikwerda}, \cite{Gustaffson}, which often result from a Fourier or von Neumann type analysis of difference schemes in the linear case. These conditions, and similar conditions for schemes to be discussed shortly, are given in the second column of Table \ref{tab4}. 

\begin{table}
\caption{\label{tab4} Stability conditions and bounds on solution size for the schemes given in Table \ref{tab1}. Schemes not possessing such conditions are indicated by the symbol $- -$ in the accompanying rows. Two types of bounds are given: in the third column, general bounds on the size of the solution, which hold under any conservative boundary conditions, and in the fourth column, better bounds available, when one of the ends of the string is fixed.} 
\begin{center}\begin{scriptsize}
\begin{tabular}{l|l|l|l}
&Stab. conditions& General bounds & Bounds under fixed conditions\\\hline\hline
${\bf s}^{(a)}$ &$- -$&$- -$&$- -$ \\\hline
${\bf s}^{(b)}$ &$- -$&$- -$&$- -$  \\\hline
\raisebox{-0.1in}{${\bf s}^{(c)}$} & $\lambda\leq \sqrt{1/\alpha}$&$\|\xi^{n}\|_{{\mathcal D}}\leq \|\xi^{0}\|_{{\mathcal D}}+h_{t}n\sqrt{\frac{2{\mathcal H}_{{\bf s}^{(c)}}^0}{1-\lambda^2}}$& $\|\xi\|_{{\mathcal D}}\leq Nh_{x}\sqrt{\frac{2{\mathcal H}_{{\bf s}^{(c)}}^0}{\alpha}}+\frac{h_{t}}{2}\sqrt{\frac{2{\mathcal H}_{{\bf s}^{(c)}}^0}{1-\lambda^2}}$ \\
&$\lambda\leq 1$&$\|\etab^{n}\|_{{\mathcal D}}\leq \|\etab^{0}\|_{{\mathcal D}}+h_{t}n\sqrt{\frac{2{\mathcal H}_{{\bf s}^{(c)}}^0}{1-\alpha\lambda^2}}$&$\|\etab\|_{{\mathcal D}}\leq Nh_{x}\sqrt{\frac{2{\mathcal H}_{{\bf s}^{(c)}}^0}{\alpha}}+\frac{h_{t}}{2}\sqrt{\frac{2{\mathcal H}_{{\bf s}^{(c)}}^0}{1-\alpha\lambda^2}}$\\\hline
\raisebox{-0.1in}{${\bf s}^{(d)}$} & $\lambda\leq \sqrt{1/\alpha}$&$\|\xi^{n}\|_{{\mathcal D}}\leq \|\xi^{0}\|_{{\mathcal D}}+h_{t}n\sqrt{\frac{2{\mathcal H}_{{\bf s}^{(d)}}^0}{1-\lambda^2}}$& $\|\xi\|_{{\mathcal D}}\leq Nh_{x}\sqrt{\frac{2{\mathcal H}_{{\bf s}^{(c)}}^0}{\alpha}}+\frac{h_{t}}{2}\sqrt{\frac{2{\mathcal H}_{{\bf s}^{(d)}}^0}{1-\lambda^2}}$\\
&$\lambda\leq 1$&$\|\etab^{n}\|_{{\mathcal D}}\leq \|\etab^{0}\|_{{\mathcal D}}+h_{t}n\sqrt{\frac{2{\mathcal H}_{{\bf s}^{(d)}}^0}{1-\alpha\lambda^2}}$&$\|\etab\|_{{\mathcal D}}\leq Nh_{x}\sqrt{\frac{2{\mathcal H}_{{\bf s}^{(d)}}^0}{\alpha}}+\frac{h_{t}}{2}\sqrt{\frac{2{\mathcal H}_{{\bf s}^{(d)}}^0}{1-\alpha\lambda^2}}$\\\hline
\raisebox{-0.1in}{${\bf s}^{(e)}$} &$\lambda\leq \sqrt{1/\alpha}$&$\|\xi^{n}\|_{{\mathcal D}}\leq \|\xi^{0}\|_{{\mathcal D}}+h_{t}n\sqrt{\frac{2{\mathcal H}_{{\bf s}^{(e)}}^0}{1-(2-\alpha)\lambda^2}}$& $\|\xi\|_{{\mathcal D}}\leq Nh_{x}\sqrt{\frac{2{\mathcal H}_{{\bf s}^{(e)}}^0}{\alpha}}+\frac{h_{t}}{2}\sqrt{\frac{2{\mathcal H}_{{\bf s}^{(e)}}^0}{1-(2-\alpha)\lambda^2}}$\\
&$\lambda\leq\sqrt{\frac{1}{2-\alpha}}$&$\|\etab^{n}\|_{{\mathcal D}}\leq \|\etab^{0}\|_{{\mathcal D}}+h_{t}n\sqrt{\frac{2{\mathcal H}_{{\bf s}^{(e)}}^0}{1-\alpha\lambda^2}}$&$\|\etab\|_{{\mathcal D}}\leq Nh_{x}\sqrt{\frac{2{\mathcal H}_{{\bf s}^{(e)}}^0}{\alpha}}+\frac{h_{t}}{2}\sqrt{\frac{2{\mathcal H}_{{\bf s}^{(e)}}^0}{1-\alpha\lambda^2}}$ \\\hline 
${\bf k}^{(a)}$ &$- -$&$- -$&$- -$\\\hline
${\bf k}^{(b)}$ &$\lambda\leq \sqrt{1/\alpha}$&$\|\etab^{n}\|_{{\mathcal D}}\leq \|\etab^{0}\|_{{\mathcal D}}+h_{t}n\sqrt{\frac{2{\mathcal H}_{{\bf k}^{(b)}}^0}{1-\alpha\lambda^2}}$& $\|\etab\|_{{\mathcal D}}\leq Nh_{x}\sqrt{\frac{2{\mathcal H}_{{\bf k}^{(b)}}^0}{\alpha}}+\frac{h_{t}}{2}\sqrt{\frac{2{\mathcal H}_{{\bf k}^{(b)}}^0}{1-\alpha\lambda^2}}$\\
 \end{tabular}
\end{scriptsize}\end{center}
\end{table}

Referring to Table \ref{tab3}, the expression ${\mathcal V}_{{\bf s}^{(c)}}$ for the discrete potential energy for scheme ${\bf s}^{(c)}$ differs from ${\mathcal V}_{{\bf s}^{(d)}}$ only by a single term, which is non-negative; following steps similar to the above, one may easily derive the following inequality:
\begin{eqnarray}
\label{scegy}
{\mathcal H}_{{\bf s}^{(c)}}&\geq& \frac{\alpha}{2}\left(\|\mu_{t-}p\|^2_{{\mathcal D}^{+}}+ \|\mu_{t-}{\bf q}\|^2_{{\mathcal D}^{+}}\right)\\
\qquad&&+\frac{1-\alpha}{2}\left(\|\mu_{t-}p+\frac{1}{2}{\bf q}^{T}e_{t-}{\bf q}\|^{2}_{{\mathcal D}^{+}}+\frac{1}{4}\|\tilde{{\bf q}}^{T}e_{t-}{\bf q}\|_{{\mathcal D}^{+}}^2\right)\notag\\
\qquad&&+ \left(\frac{1}{2}-\frac{\lambda^2}{2}\right)\|\delta_{t-}\xi\|^2_{{\mathcal D}}+\left(\frac{1}{2}-\frac{\alpha\lambda^2}{2}\right)\|\delta_{t-}\etab\|^2_{{\mathcal D}}\notag
\end{eqnarray}
Thus the positivity conditions \eqref{CFLcond} derived above for scheme ${\bf s}^{(d)}$ hold for scheme ${\bf s}^{(c)}$ as well.

 The analysis for scheme ${\bf s}^{(e)}$ is similar to that performed above for schemes ${\bf s}^{(d)}$, and ${\bf s}^{(c)}$, but the positivity conditions are slightly different. One may derive the following inequality:
\begin{eqnarray}
\label{seegy}
{\mathcal H}_{{\bf s}^{(e)}}&\geq& \frac{\alpha}{2}\left(\|\mu_{t-}p\|^2_{{\mathcal D}^{+}}+ \|\mu_{t-}{\bf q}\|^2_{{\mathcal D}^{+}}\right)\\
\qquad&&+\frac{1-\alpha}{2}\mu_{t-}\left(\|p+\frac{1}{2}{\bf q}^{T}{\bf q}\|^{2}_{{\mathcal D}^{+}}\right)\notag\\
\qquad&&+ \left(\frac{1}{2}-\frac{(2-\alpha)\lambda^2}{2}\right)\|\delta_{t-}\xi\|^2_{{\mathcal D}}+\left(\frac{1}{2}-\frac{\alpha\lambda^2}{2}\right)\|\delta_{t-}\etab\|^2_{{\mathcal D}}\notag
\end{eqnarray}
For positivity, condition \eqref{CFLcond2} holds as before, but condition \eqref{CFLcond1} must be modified to
\begin{equation}
\label{CFLcond3}
\lambda\leq\sqrt{\frac{1}{2-\alpha}}
\end{equation}

For schemes ${\bf s}^{(c)}$, ${\bf s}^{(d)}$, and ${\bf s}^{(e)}$, the analysis above is simple, because the contributions of the nonlinearity to the expressions for the discrete conserved energy (i.e., the terms which are not simply quadratic in the expressions for ${\mathcal V}_{{\bf s}^{(c)}}$, ${\mathcal V}_{{\bf s}^{(d)}}$ and ${\mathcal V}_{{\bf s}^{(e)}}$) are themselves non-negative. Thus the determination of positivity conditions reduces, essentially, to analysis of a linear problem. For scheme ${\bf s}^{(b)}$, however, this is not the case; the term in the expression for ${\mathcal V}_{{\bf s}^{(b)}}$ resulting from the nonlinearity is not necessarily positive, and indeed, can be negative and unbounded (the choices $p=0$ and $e_{t-}p = -e_{t-}{\bf q}^{T}{\bf q}$ illustrate this point simply). No global positivity condition is available for scheme ${\bf s}^{(b)}$. 

For the energy-conserving scheme ${\bf k}^{(b)}$, similar positivity conditions are also immediate. Beginning from the expression for ${\mathcal V}_{{\bf k}^{(b)}}$, one may write, employing identity \eqref{shiftid} and inequality \eqref{diffbound}, 
\begin{eqnarray*}
{\mathcal V}_{{\bf k}^{(b)}} &=& \frac{\alpha}{2}\|\mu_{t-}{\bf q}\|_{{\mathcal D}^{+}}^2-\frac{\alpha h_{t}^2}{8}\|\delta_{t-}{\bf q}\|_{{\mathcal D}^{+}}^2+\frac{1}{8}\langle {\bf q}, e_{t-}{\bf q}\rangle_{{\mathcal D}^{+}}^2\\
&\geq& \frac{\alpha}{2}\|\mu_{t-}{\bf q}\|_{{\mathcal D}^{+}}^2-\frac{\alpha\lambda^2}{2}\|\delta_{t-}\etab\|_{{\mathcal D}^{+}}^2+\frac{1}{8}\langle {\bf q}, e_{t-}{\bf q}\rangle_{{\mathcal D}^{+}}^2
\end{eqnarray*}
implying
\begin{equation}
\label{kbegy}
{\mathcal H}_{{\bf k}^{(b)}} = {\mathcal T}_{{\bf k}^{(b)}}+ {\mathcal V}_{{\bf k}^{(b)}} \geq \frac{\alpha}{2}\|\mu_{t-}{\bf q}\|_{{\mathcal D}^{+}}^2+\left(\frac{1}{2}-\frac{\alpha\lambda^2}{2}\right)\|\delta_{t-}\etab\|_{{\mathcal D}^{+}}^2+\frac{1}{8}\langle {\bf q}, e_{t-}{\bf q}\rangle_{{\mathcal D}^{+}}^2
\end{equation}
The total energy is then positive under condition \eqref{CFLcond1}. 

\subsubsection{Bounds on Solution Size}

In this section, only schemes with a conditionally positive discrete conserved energy are considered, namely schemes ${\bf s}^{(c)}$, ${\bf s}^{(d)}$, ${\bf s}^{(e)}$ and ${\bf k}^{(b)}$, under conservative boundary conditions.  

Consider first the inequalities \eqref{scegy} and \eqref{sdegy}, which give lower bounds on the conserved energy for schemes ${\bf s}^{(c)}$ and ${\bf s}^{(d)}$. Under stability conditions \eqref{CFLcond}, it is then clear that, for either scheme,
\begin{equation}
\label{genbound}
\|\delta_{t-}\xi\|_{{\mathcal D}} \leq \sqrt{\frac{2{\mathcal H}_{{\bf s}^{(\bullet)}}^0}{1-\lambda^2}}\qquad\|\delta_{t-}\etab\|_{{\mathcal D}} \leq \sqrt{\frac{2{\mathcal H}_{{\bf s}^{(\bullet)}}^0}{1-\alpha\lambda^2}}
\end{equation}
For scheme ${\bf s}^{(e)}$, the bounds, under stability conditions \eqref{CFLcond2} and \eqref{CFLcond3}, are slightly different; now, one has
\begin{equation}
\label{genbounde}
\|\delta_{t-}\xi\|_{{\mathcal D}} \leq \sqrt{\frac{2{\mathcal H}_{{\bf s}^{(e)}}^0}{1-(2-\alpha)\lambda^2}}\qquad\|\delta_{t-}\etab\|_{{\mathcal D}} \leq \sqrt{\frac{2{\mathcal H}_{{\bf s}^{(e)}}^0}{1-\alpha\lambda^2}}
\end{equation}
For scheme ${\bf k}^{(b)}$, the bound on \etab\, follows immediately as
\begin{equation*}
\|\delta_{t-}\etab\|_{{\mathcal D}} \leq \sqrt{\frac{2{\mathcal H}_{{\bf k}^{(b)}}^0}{1-\alpha\lambda^2}}
\end{equation*}

Any of these bounds on the norm of a time difference of a grid function may be converted to a general bound on the norm of the grid function itself, through an application of inequality \eqref{tdbound}. For example, for the bounds given above on $\delta_{t-}\xi$ in \eqref{genbound} for schemes ${\bf s}^{(c)}$ or ${\bf s}^{(d)}$, it then follows that
\begin{equation*}
\|\xi^{n}\|_{{\mathcal D}}\leq \|\xi^{0}\|_{{\mathcal D}}+h_{t}n\sqrt{\frac{2{\mathcal H}_{{\bf s}^{(\bullet)}}^0}{1-\lambda^2}}
\end{equation*}
Thus growth of the longitudinal displacement is at most linear in time, for any conservative boundary conditions. A similar bound can be found for the transverse displacement, using the second of bounds \eqref{genbound}; bounds for the schemes ${\bf s}^{(e)}$ and ${\bf k}^{(b)}$, arrived at in a nearly identical way, are given in Table \ref{tab4}. 

Under fixed boundary conditions, tighter bounds may be obtained. For schemes ${\bf s}^{(c)}$, ${\bf s}^{(d)}$, and ${\bf s}^{(e)}$, the following bounds on the quantities $p$ and ${\bf q}$ follow from the appropriate stability conditions (\eqref{CFLcond} for ${\bf s}^{(c)}$ and ${\bf s}^{(d)}$, and \eqref{CFLcond2} and \eqref{CFLcond3} for ${\bf s}^{(e)}$) and inequalities \eqref{scegy}, \eqref{sdegy} and \eqref{seegy}:
\begin{equation}
\label{fixedbound}
\|\mu_{t-}p\|_{{\mathcal D}^{+}} \leq \sqrt{\frac{2{\mathcal H}_{{\bf s}^{(\bullet)}}^0}{\alpha}}\qquad\|\mu_{t-}{\bf q}\|_{{\mathcal D}^{+}} \leq \sqrt{\frac{2{\mathcal H}_{{\bf s}^{(\bullet)}}^0}{\alpha}}
\end{equation}
For scheme ${\bf k}^{(b)}$, the bound, from stability condition \eqref{CFLcond2} and inequality \eqref{kbegy}, is simply 
\begin{equation*}
\|\mu_{t-}{\bf q}\|_{{\mathcal D}^{+}} \leq \sqrt{\frac{2{\mathcal H}_{{\bf s}^{(\bullet)}}^0}{\alpha}}
\end{equation*}

For any of the above bounds on $\|\mu_{t-}p\|_{{\mathcal D}^{+}}$ or $\|\mu_{t-}{\bf q}\|_{{\mathcal D}^{+}}$, bounds on $\|\xi\|_{{\mathcal D}}$ or $\|\etab\|_{{\mathcal D}}$ may be obtained in the following way. Considering, for example, \etab\,, one may immediately write, using identity \eqref{mudid} and the triangle inequality \eqref{trieq}, 
\begin{equation*}
\etab = \mu_{t-}\etab+\frac{h_{t}}{2}\delta_{t-}\etab\qquad\Longrightarrow\qquad\|\etab\|_{{\mathcal D}} \leq \|\mu_{t-}\etab\|_{{\mathcal D}}+\frac{h_{t}}{2}\|\delta_{t-}\etab\|_{{\mathcal D}}
\end{equation*}
and, furthermore, under scheme ${\bf s}^{(c)}$ for example, 
\begin{equation*}
\|\etab\|_{{\mathcal D}} \leq Nh_{x}\|\mu_{t-}{\bf q}\|_{{\mathcal D}^{+}}+\frac{h_{t}}{2}\sqrt{\frac{2{\mathcal H}_{{\bf s}^{(c)}}^0}{1-\alpha\lambda^2}}\leq Nh_{x}\sqrt{\frac{2{\mathcal H}_{{\bf s}^{(c)}}^0}{\alpha}}+\frac{h_{t}}{2}\sqrt{\frac{2{\mathcal H}_{{\bf s}^{(c)}}^0}{1-\alpha\lambda^2}}
\end{equation*}
where in the first inequality above, the inequality \eqref{spatubound} and the second of the bounds \eqref{genbound} have been employed, and in the second, the second of the bounds \eqref{fixedbound} has been used. 

This bound, and similar forms for the other schemes of interest in this section are given in Table \ref{tab4}. 
\subsection{Implementation Details}
\label{implementationsec}
It is useful to write the schemes defined in Table \ref{tab1} in forms suitable for computer implementation. In order to avoid presenting a multiplicity of slightly different cases, it is assumed in this section that boundary conditions are of the fixed type (i.e., \eqref{bcdef1fd} for system {\bf s}) at both ends of the string, which are the conditions of interest in many cases. In this case, the values of $\etab_{i}^{n}$ and $\xi_{i}^{n}$ at the endpoints of the spatial domain (i.e., at $i=0$ and $i=N$) are fixed at zero, and need not be considered. It is helpful to introduce the column vectors ${\bf u}^{n}$, ${\bf v}_{(1)}^n$ and ${\bf v}_{(2)}^n$ defined by 
\begin{equation*}
{\bf u}^{n} = [\xi_{1}^{n},\hdots,\xi_{N-1}^{n}]^{T}\qquad{\bf v}_{(1)}^n = [\eta_{(1),1}^n, \hdots, \eta_{(1),N-1}^n]^{T}\qquad{\bf v}_{(2)}^n = [\eta_{(2),1}^{n}, \hdots, \eta_{(2),N-1}^{n}]^{T}
\end{equation*}
and the state vector ${\bf w}^{n}$ containing all the displacements at time step $n$, defined by
\begin{equation*}
{\bf w}^{n} = [({\bf u}^{n})^{T}, ({\bf v}_{(1)}^{n})^{T}, ({\bf v}_{(2)}^{n})^{T}]^{T}
\end{equation*}
Schemes ${\bf s}^{(a)}$, ${\bf s}^{(b)}$, ${\bf s}^{(c)}$, and ${\bf s}^{(d)}$ may all be written in the form of two-step matrix recursions of the form 
\begin{equation}
\label{updateform}
{\bf w}^{n+1}=({\bf A}_{{\bf s}^{(\bullet)}}^{n})^{-1}{\bf B}_{{\bf s}^{(\bullet)}}^{n}{\bf w}^{n}-{\bf w}^{n-1}
\end{equation}
where ${\bf A}_{{\bf s}^{(\bullet)}}^{n}$ and ${\bf B}_{{\bf s}^{(\bullet)}}^{n}$ are square matrices which depend on the values of the state at time step $n$, i.e., ${\bf w}^{n}$. The fact that ${\bf A}_{{\bf s}^{(\bullet)}}^{n}$ must be inverted (or, rather, a linear system solved) reflects the implicit nature of the schemes. For system ${\bf s}^{(a)}$, the matrix ${\bf A}_{{\bf s}^{(\bullet)}}^{n}$ is diagonal, and the algorithm is explicit. The complete forms of the state update matrices are given in Table \ref{tab5}, where the scaled differentiation matrix ${\bf D}$ is defined as
\begin{equation}
\label{ddef}
{\bf D} = \lambda\begin{bmatrix}
1 & & &  \\
-1 & 1 &  \\
 & \ddots & \ddots & \\
 &  & -1 & 1\\
 &  & & -1\\
 \end{bmatrix}
 \end{equation}
 ${\bf D}$ is an $N$ by $N-1$ matrix. In addition, it is useful to define the diagonal $N$ by $N$ matrices ${\bf P}$, ${\bf Q}_{(1)}$ and ${\bf Q}_{(2)}$ by
\begin{equation*}
{\bf P} = \frac{1}{h_{t}}\mbox{{\rm diag}}({\bf D}{\bf u}^{n})\qquad{\bf Q}_{(1)} = \frac{1}{h_{t}}\mbox{{\rm diag}}({\bf D}{\bf v}_{(1)}^{n}\qquad{\bf Q}_{(2)} = \frac{1}{h_{t}}\mbox{{\rm diag}}({\bf D}{\bf v}_{(2)}^{n})
\end{equation*}
In Table \ref{tab5}, ${\bf I}$ refers to the $N-1$ by $N-1$ identity matrix, and the constant $\beta$ is defined as
\begin{equation*}
\beta = \frac{\alpha-1}{2}
\end{equation*} 

\begin{table}
\caption{\label{tab5} State update matrices for schemes ${\bf s}^{(\bullet)}$, under fixed boundary conditions. The block matrices ${\bf A}_{{\bf s}^{(\bullet)}}$ and ${\bf B}_{{\bf s}^{(\bullet)}}$, given below, are used in the update of Eq. \eqref{updateform}. The symbol $\cdot$ indicates a an $N-1$ by $N-1$ zero matrix. System ${\bf s}^{(e)}$ does not have a state update form. } 
\begin{center}\begin{scriptsize}
\begin{tabular}{l|l}
& \\\hline\hline
${\bf s}^{(a)}$ & ${\bf A}_{{\bf s}^{(a)}} = \begin{bmatrix}{\bf I}&\cdot&\cdot\\\cdot&{\bf I}&\cdot\\\cdot&\cdot&{\bf I}\\ \end{bmatrix}$ \\
& ${\bf B}_{{\bf s}^{(a)}} = \begin{bmatrix} 2{\bf I}-{\bf D}^{T}{\bf D}& \beta{\bf D}^{T}{\bf Q}^{(1)}{\bf D}& \beta{\bf D}^{T}{\bf Q}^{(2)}{\bf D}\\ \beta{\bf D}^{T}{\bf Q}^{(1)}{\bf D} & 2{\bf I}-\alpha{\bf D}^{T}{\bf D}+\beta{\bf D}^{T}\left({\bf Q}_{(1)}^2+{\bf Q}_{(2)}^2+{\bf P}\right){\bf D} & \cdot\\\beta{\bf D}^{T}{\bf Q}^{(2)}{\bf D} & \cdot & 2{\bf I}-\alpha{\bf D}^{T}{\bf D}+\beta{\bf D}^{T}\left({\bf Q}_{(1)}^2+{\bf Q}_{(2)}^2+{\bf P}\right){\bf D} \end{bmatrix}$\\\hline
${\bf s}^{(b)}$ & ${\bf A}_{{\bf s}^{(b)}} = \begin{bmatrix}{\bf I}&\cdot&\cdot\\\cdot&{\bf I}-\frac{\beta}{2}{\bf D}^{T}\left({\bf Q}_{(1)}^2+{\bf Q}_{(2)}^2+2{\bf P}\right){\bf D}&\cdot\\\cdot&\cdot&{\bf I}-\frac{\beta}{2}{\bf D}^{T}\left({\bf Q}_{(1)}^2+{\bf Q}_{(2)}^2+2{\bf P}\right){\bf D}\\ \end{bmatrix}$\\
& ${\bf B}_{{\bf s}^{(b)}} = \begin{bmatrix} 2{\bf I}-{\bf D}^{T}{\bf D}& \beta{\bf D}^{T}{\bf Q}_{(1)}{\bf D}& \beta{\bf D}^{T}{\bf Q}_{(2)}{\bf D}\\ \cdot & 2{\bf I}-\alpha{\bf D}^{T}{\bf D} & \cdot\\\cdot & \cdot &  2{\bf I}-\alpha{\bf D}^{T}{\bf D} \end{bmatrix}$\\\hline
${\bf s}^{(c)}$ & ${\bf A}_{{\bf s}^{(c)}} = \begin{bmatrix}{\bf I}&-\frac{\beta}{2}{\bf D}^{T}{\bf Q}_{(1)}{\bf D}&-\frac{\beta}{2}{\bf D}^{T}{\bf Q}_{(2)}{\bf D}\\-\frac{\beta}{2}{\bf D}^{T}{\bf Q}_{(1)}{\bf D}&{\bf I}-\frac{\beta}{2}{\bf D}^{T}\left({\bf Q}_{(1)}^2+{\bf Q}_{(2)}^2\right){\bf D}&\cdot\\-\frac{\beta}{2}{\bf D}^{T}{\bf Q}_{(2)}{\bf D}&\cdot&{\bf I}-\frac{\beta}{2}{\bf D}^{T}\left({\bf Q}_{(1)}^2+{\bf Q}_{(2)}^2\right){\bf D}\\ \end{bmatrix}$\\
& ${\bf B}_{{\bf s}^{(c)}} = \begin{bmatrix} 2{\bf I}-{\bf D}^{T}{\bf D}& \cdot& \cdot\\ \cdot & 2{\bf I}-\alpha{\bf D}^{T}{\bf D}+\beta{\bf D}^{T}{\bf P}{\bf D}  & \cdot\\\cdot & \cdot &  2{\bf I}-\alpha{\bf D}^{T}{\bf D}+\beta{\bf D}^{T}{\bf P}{\bf D} \end{bmatrix}$\\\hline
${\bf s}^{(d)}$ & ${\bf A}_{{\bf s}^{(d)}} = \begin{bmatrix}{\bf I}&-\frac{\beta}{2}{\bf D}^{T}{\bf Q}_{(1)}{\bf D}&-\frac{\beta}{2}{\bf D}^{T}{\bf Q}_{(2)}{\bf D}\\-\frac{\beta}{2}{\bf D}^{T}{\bf Q}_{(1)}{\bf D}&{\bf I}-\frac{\beta}{2}{\bf D}^{T}\left({\bf Q}_{(1)}^2\right){\bf D}&-\frac{\beta}{2}{\bf D}^{T}\left({\bf Q}_{(1)}{\bf Q}_{(2)}\right){\bf D}\\-\frac{\beta}{2}{\bf D}^{T}{\bf Q}_{(2)}{\bf D}&-\frac{\beta}{2}{\bf D}^{T}\left({\bf Q}_{(1)}{\bf Q}_{(2)}\right){\bf D}&{\bf I}-\frac{\beta}{2}{\bf D}^{T}\left({\bf Q}_{(2)}^2\right){\bf D}\\ \end{bmatrix}$\\
& ${\bf B}_{{\bf s}^{(d)}} = \begin{bmatrix} 2{\bf I}-{\bf D}^{T}{\bf D}& \cdot& \cdot\\\cdot  & 2{\bf I}-\alpha{\bf D}^{T}{\bf D}+\beta{\bf D}^{T}{\bf P}{\bf D}  & \cdot\\\cdot & \cdot &  2{\bf I}-\alpha{\bf D}^{T}{\bf D}+\beta{\bf D}^{T}{\bf P}{\bf D} \end{bmatrix}$\\\hline
${\bf s}^{(e)}$ & \mbox{{\rm No matrix update form}}\\\hline
\end{tabular}
\end{scriptsize}\end{center}
\end{table}

If a scheme may be written in a state update form as per Eq. \eqref{updateform}, one has an immediate proof of existence and uniqueness of solutions, and this is indeed the case for schemes ${\bf s}^{(a)}$, ${\bf s}^{(b)}$, ${\bf s}^{(c)}$, and ${\bf s}^{(d)}$. In all these cases, although the difference scheme is, as a whole, strongly nonlinear and implicit, the coupling among the state variables at the current time step is linear in character. It is important to note, however, that this property is independent of the existence of discrete conservation laws---for scheme ${\bf s}^{(e)}$, which is energy conserving, there is no state update form. Thus a proof of conditions for existence and uniqueness for solutions to ${\bf s}^{(e)}$ will be considerably more difficult, if even possible to obtain. It remains true, however, that if a solution does exist, it will be stable under condition \eqref{CFLcond}. 

For scheme ${\bf k}^{(a)}$, which is explicit, the implementation is immediate. For scheme ${\bf k}^{(b)}$, however, the situation is slightly different. Though it might appear, from the form of ${\mathcal G}_{{\bf k}^{(b)}}$, that the scheme would be implicit, it is in fact possible to write it in an explicit form in the following way. For the quantity ${\mathcal G}_{{\bf k}^{(b)}}$, one may proceed as follows:
\begin{eqnarray*}
{\mathcal G}_{{\bf k}^{(b)}} &=& 1+\frac{1}{2\alpha}\mu_{t+}\langle{\bf q}, e_{t-}{\bf q}\rangle_{D^{+}}\\
&=& 1+\frac{1}{2\alpha}\langle{\bf q}, \mu_{to}{\bf q}\rangle_{D^{+}}\\
&=& 1+\frac{1}{2\alpha}\|{\bf q}\|_{D^{+}}^2+\frac{h_{t}^2}{4\alpha}\langle{\bf q},\delta_{t+}\delta_{t-}{\bf q}\rangle_{D^{+}}\\
&=& 1+\frac{1}{2\alpha}\|{\bf q}\|_{D^{+}}^2+\frac{h_{t}^2}{4}{\mathcal G}_{{\bf k}^{(b)}}\langle{\bf q},\delta_{x+}\delta_{x-}{\bf q}\rangle_{D^{+}}
\end{eqnarray*}
where in the second, third and fourth equalities, identity \eqref{shiftid2}, the definition of system ${\bf k}^{(b)}$ from Table \ref{tab1}, and identity \eqref{shitid2} have been used, respectively. 
and thus
\begin{equation*}
{\mathcal G}_{{\bf k}^{(b)}} = \frac{1+\frac{1}{2\alpha}\|{\bf q}\|_{D^{+}}}{1-\frac{h_{t}^2}{4}\langle{\bf q},\delta_{x+}\delta_{x-}{\bf q}\rangle_{D^{+}}}
\end{equation*}
which renders scheme ${\bf k}^{(b)}$ fully explicit. If the boundary conditions are fixed or free, one may go further and write
\begin{equation*}
{\mathcal G}_{{\bf k}^{(b)}} = \frac{1+\frac{1}{2\alpha}\|{\bf q}\|_{D^{+}}^2}{1+\frac{h_{t}^2}{4}\|\delta_{x+}{\bf q}\|_{[1,N-1]}^2}
\end{equation*}

\section{Numerical Examples}
\label{numsec}
\subsection{Energy and Angular Momentum Conservation}
\label{numconssec}
\begin{figure}[ht]
\centerline{\includegraphics[scale=0.75,clip,trim=44mm 96mm 48mm 96mm]{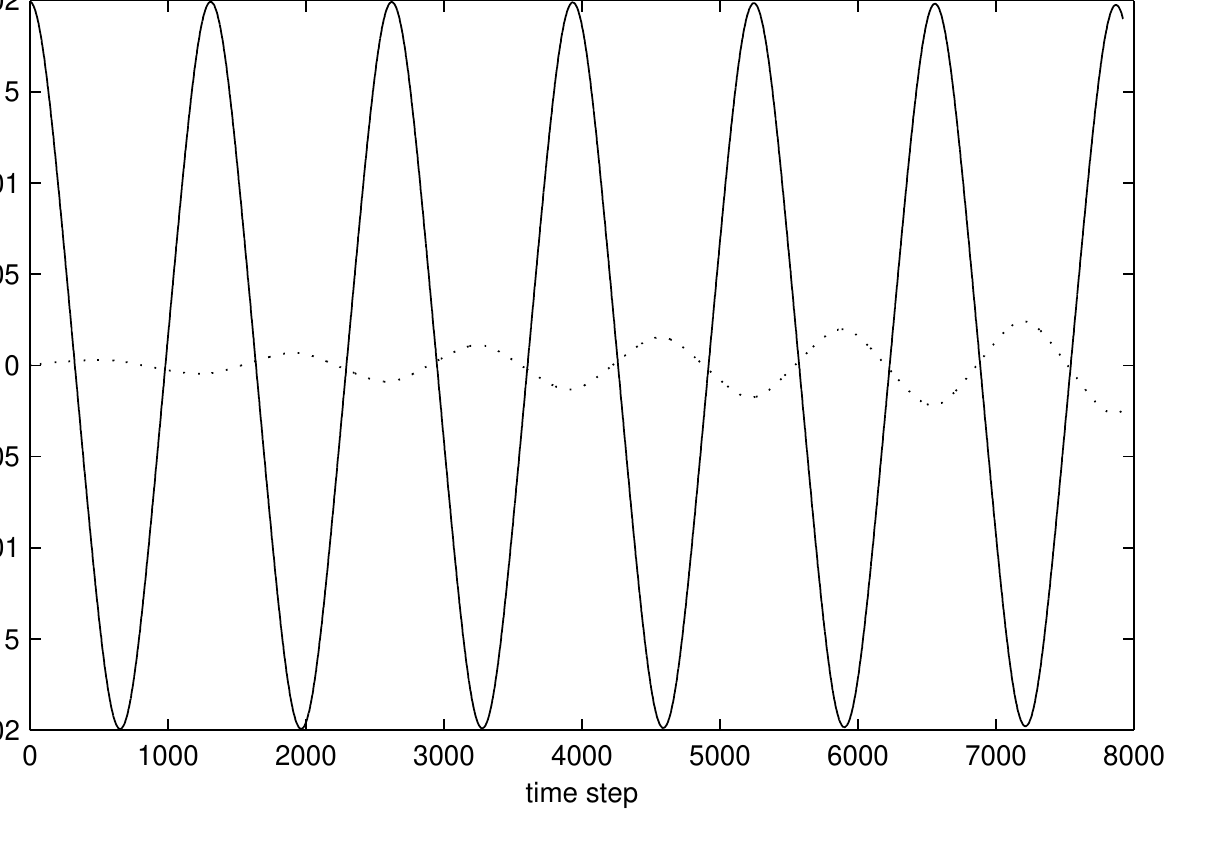}\put(-290,100){{\small $\etab$}}}
\caption{Transverse displacement at center of string (with $\alpha = 2\times 10^{-4}$), with fixed boundary conditions, plotted against time step, for schemes ${\bf s}^{(a)}$, under initial conditions given in Eqs. \eqref{ex1}, with $\gamma_{1}=0.02$ and $\gamma_{2}=2\times 10^{-5}$. The time step was chosen as $h_{t} = 1/20$. Displacement component $\eta_{(1)}$ is shown as a solid line, and $\eta_{(2)}$ as a dotted line.}  
\label{fig1}
\end{figure}

As a basic illustration of the conservative properties of the schemes presented here, consider the case of a steel string string, of parameters given in the caption to Fig. \ref{fig1}. The initial conditions of the string are 
\begin{equation}
\label{ex1}
\eta_{(1)}(x,0) = \gamma_{1}\sin(\pi x)\qquad \dot{\eta}_{(2)}(x,0) = \gamma_{2}\sin(\pi x)
\end{equation}
and are discretized as
\begin{equation}
\label{ex1a}
\eta_{(1),i}^0 = \eta_{(1),i}^1=\gamma_{1}\sin(\pi ih_{x})\qquad\eta_{(2),i}^0 = 0\qquad\eta_{(2),i}^1=h_{t}\gamma_{2}\sin(\pi ih_{x})
\end{equation}
In other words, the string is subjected to an initial displacement in the form of the first linear modal configuration in the polarization corresponding to $\eta_{(1)}$, and an initial velocity of a similar form in the other polarization. The initial longitudinal displacement and velocity is assumed to be zero. Shown in Fig. \ref{fig1} are the results of simulations using scheme ${\bf s}^{(a)}$, where the transverse string displacements $\eta_{(1)}$ and $\eta_{(2)}$ at the string center are plotted as a function of time in each case. In this, and all examples in this section, the parameter $\lambda$ is chosen as close to the bounds \eqref{CFLcond2} as possible.  

Tables \ref{tab6} and \ref{tab7} illustrate the numerical conservation properties of those algorithms which possess discrete conserved angular momentum and energy, respectively. The fluctuations observed in some cases (typically in the 12th place) are due to numerical round-off error. 

\begin{table}
\caption{\label{tab6} Conserved angular momentum (nondimensional, $\times 10^{-7}$) for various algorithms for the string described in the caption of Fig. \ref{fig1}, against time step $n$.} 
\begin{center}\begin{scriptsize}
\begin{tabular}{l|l|l|l|l|l}
n &${\mathcal A}_{{\bf s}^{(a)}}$ & ${\mathcal A}_{{\bf s}^{(b)}}$& ${\mathcal A}_{{\bf s}^{(d)}}$&${\mathcal A}_{{\bf k}^{(a)}}$ & ${\mathcal A}_{{\bf k}^{(b)}}$\\\hline\hline
1 &2.0000000000000&1.9999950753362&2.0000000000000&2.0000000000000&2.0000000000000\\\hline
2 &2.0000000000000&1.9999950753362&2.0000000000000&2.0000000000000&2.0000000000000\\\hline
3 &2.0000000000000&1.9999950753362&2.0000000000000&2.0000000000000&2.0000000000000\\\hline
4 &2.0000000000000&1.9999950753362&2.0000000000000&2.0000000000000&2.0000000000000\\\hline
5 &2.0000000000000&1.9999950753362&2.0000000000000&2.0000000000000&2.0000000000000\\\hline
100 &2.0000000000000&1.9999950753362&2.0000000000007&2.0000000000000&2.0000000000001\\
\end{tabular}
\end{scriptsize}\end{center}
\end{table}

\begin{table}
\caption{\label{tab7} Conserved energy (nondimensional, $\times 10^{-7}$) for various algorithms for the string described in the caption of Fig. \ref{fig1}, against time step $n$.} 
\begin{center}\begin{scriptsize}
\begin{tabular}{l|l|l|l|l}
n &${\mathcal H}_{{\bf s}^{(b)}}$ & ${\mathcal H}_{{\bf s}^{(c)}}$& ${\mathcal H}_{{\bf s}^{(d)}}$&${\mathcal H}_{{\bf k}^{(b)}}$\\\hline\hline
1 &9.245104334637&9.245104334637&9.245104316451&6.821328138420\\\hline
2 &9.245104334637&9.245104334637&9.245104316451&6.821328138420\\\hline
3 &9.245104334637&9.245104334637&9.245104316451&6.821328138420\\\hline
4 &9.245104334637&9.245104334637&9.245104316451&6.821328138420\\\hline
5 &9.245104334637&9.245104334637&9.245104316451&6.821328138420\\\hline
100 &9.245104334635&9.245104334635&9.245104316452&6.821328138419\\\hline\end{tabular}
\end{scriptsize}\end{center}
\end{table}

\begin{figure}[ht]
\centerline{\includegraphics[scale=0.70,clip,trim=14mm 120mm 8mm 76mm]{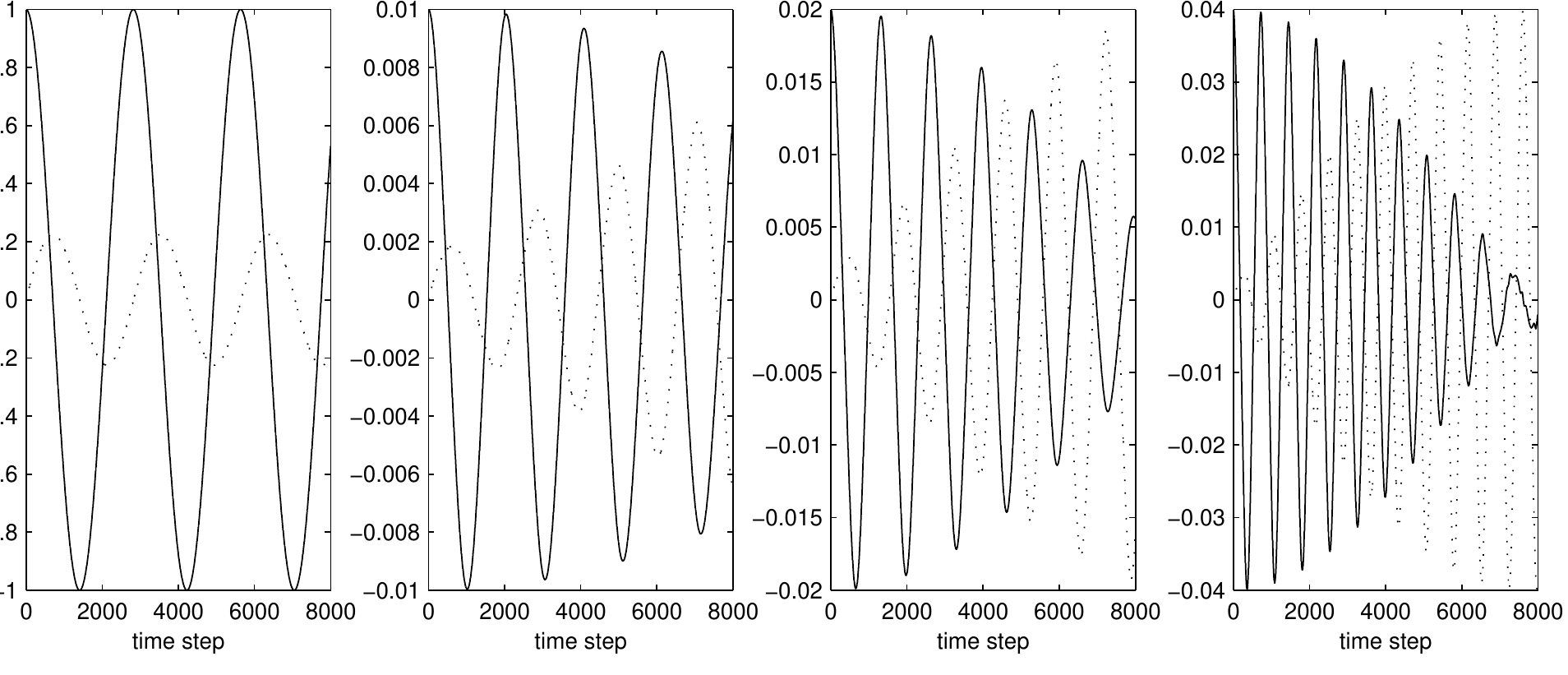}\multiput(-400,95)(99,0){4}{{\scriptsize $\etab$}}
\put(-355,180){{\scriptsize $\gamma_{1}=0.001$}}\put(-355,170){{\scriptsize $\gamma_{2}=0.00001$}}\put(-255,180){{\scriptsize $\gamma_{1}=0.01$}}\put(-255,170){{\scriptsize $\gamma_{2}=0.0001$}}\put(-155,180){{\scriptsize $\gamma_{1}=0.02$}}\put(-155,170){{\scriptsize $\gamma_{2}=0.0002$}}\put(-55,180){{\scriptsize $\gamma_{1}=0.04$}}\put(-55,170){{\scriptsize $\gamma_{2}=0.0004$}}}
\caption{Transverse displacement at center of string (with $\alpha = 2\times 10^{-4}$), with fixed boundary conditions, plotted against time step, for scheme ${\bf s}^{(d)}$, with initial conditions \eqref{ex1}, under different values of $\gamma_{1}$ and $\gamma_{2}$, as given in the panels above. The time step was chosen as $h_{t} = 1/20$. Displacement component $\eta_{(1)}$ is shown as a solid line, and $\eta_{(2)}$ as a dotted line.}  
\label{fig2}
\end{figure}

As an example of the effect of the nonlinearity, consider a representative scheme, such as, e.g., ${\bf s}^{(d)}$, which is both angular momentum and energy conserving, applied to the string of parameters as given in the caption to Fig. \ref{fig1}, and again with initial conditions as given by Eqs. \eqref{ex1}.  In Fig. \ref{fig2} are shown transverse displacements at the string center, for different values of $\gamma_{1}$ and $\gamma_{2}$. Notice in particular the change in the oscillating frequency. Values of the conserved quantities ${\mathcal A}_{{\bf s}^{(d)}}$ and ${\mathcal H}_{{\bf s}^{(d)}}$ are given in tables \ref{tab8} and \ref{tab9}, respectively. 

\begin{table}
\caption{\label{tab8} Conserved angular momentum ${\mathcal A}_{{\bf s}^{(d)}}$ (nondimensional) for scheme ${\bf s}^{(d)}$, against time step $n$, for the string of parameters defined in the caption to Fig. \ref{fig1}, under initial conditions \eqref{ex1} of increasing magnitude. Values of $\gamma_{1}$ and $\gamma_{2}$ are given in the table.}
\begin{center}\begin{scriptsize}
\begin{tabular}{l|l|l|l|l}
n &$\gamma_{1} = 0.001\,\,\gamma_{2}=0.00001$ & $\gamma_{1} = 0.01\,\,\gamma_{2}=0.0001$& $\gamma_{1} = 0.02\,\,\gamma_{2}=0.0002$&$\gamma_{1} = 0.04\,\,\gamma_{2}=0.0004$ \\\hline\hline
1 & 5.000000000000$\times 10^{-9}$ & 5.000000000000$\times 10^{-7}$ & 2.000000000000$\times 10^{-6}$ & 8.000000000000$\times 10^{-6}$\\\hline
2 & 5.000000000000$\times 10^{-9}$ & 5.000000000000$\times 10^{-7}$ & 2.000000000000$\times 10^{-6}$ & 8.000000000000$\times 10^{-6}$\\\hline
3 & 5.000000000000$\times 10^{-9}$ & 5.000000000000$\times 10^{-7}$ & 2.000000000000$\times 10^{-6}$ & 8.000000000000$\times 10^{-6}$\\\hline
4 & 5.000000000000$\times 10^{-9}$ & 5.000000000000$\times 10^{-7}$ & 2.000000000000$\times 10^{-6}$ & 8.000000000000$\times 10^{-6}$\\\hline
5 & 5.000000000000$\times 10^{-9}$ & 5.000000000000$\times 10^{-7}$ & 2.000000000000$\times 10^{-6}$ & 8.000000000000$\times 10^{-6}$\\\hline
100 & 5.000000000002$\times 10^{-9}$ & 5.000000000002$\times 10^{-7}$ & 1.999999999999$\times 10^{-6}$ & 8.000000000000$\times 10^{-6}$\\\hline
\end{tabular}
\end{scriptsize}\end{center}
\end{table}

\begin{table}
\caption{\label{tab9} Conserved energy ${\mathcal H}_{{\bf s}^{(d)}}$ (nondimensional) for scheme ${\bf s}^{(d)}$, against time step $n$, for the string of parameters defined in the caption to Fig. \ref{fig1}, under initial conditions \eqref{ex1} of increasing magnitude. Values of $\gamma_{1}$ and $\gamma_{2}$ are given in the table.}
\begin{center}\begin{scriptsize}
\begin{tabular}{l|l|l|l|l}
n &$\gamma_{1} = 0.001\,\,\gamma_{2}=0.00001$ & $\gamma_{1} = 0.01\,\,\gamma_{2}=0.0001$& $\gamma_{1} = 0.02\,\,\gamma_{2}=0.0002$&$\gamma_{1} = 0.04\,\,\gamma_{2}=0.0004$ \\\hline\hline
1 & 5.22012775452$\times 10^{-10}$ & 9.72106301924$\times 10^{-8}$ & 9.34410431645$\times 10^{-7}$ & 1.24667283005$\times 10^{-5}$ \\\hline
2 & 5.22012775452$\times 10^{-10}$ & 9.72106301924$\times 10^{-8}$ & 9.34410431645$\times 10^{-7}$ & 1.24667283005$\times 10^{-5}$ \\\hline
3 & 5.22012775452$\times 10^{-10}$ & 9.72106301924$\times 10^{-8}$ & 9.34410431645$\times 10^{-7}$ & 1.24667283005$\times 10^{-5}$ \\\hline
4 & 5.22012775452$\times 10^{-10}$ & 9.72106301924$\times 10^{-8}$ & 9.34410431645$\times 10^{-7}$ & 1.24667283005$\times 10^{-5}$ \\\hline
5 & 5.22012775452$\times 10^{-10}$ & 9.72106301924$\times 10^{-8}$ & 9.34410431645$\times 10^{-7}$ & 1.24667283005$\times 10^{-5}$ \\\hline
100 & 5.22012775452$\times 10^{-10}$ & 9.72106301925$\times 10^{-8}$ & 9.34410431645$\times 10^{-7}$ & 1.24667283005$\times 10^{-5}$ \\\hline
\end{tabular}
\end{scriptsize}\end{center}
\end{table}

\subsection{Instability of Planar Motion}
\label{planarsec}

The phenomenon of the instability of purely planar motion is inherent to nonlinear strings, and has been examined in depth by various authors \cite{Dickey80}, \cite{Gough}, \cite{Oreilly}, \cite{Rowland}. (By ``instability," one refers here to the tendency for motion which is confined to a single plane to be transferred to the perpendicular polarization, and not instability in the sense of explosive growth of the solution.) 

\begin{figure}[h!!!]
\centerline{\includegraphics[scale=0.85,clip,trim=0mm 70mm 23mm 40mm]{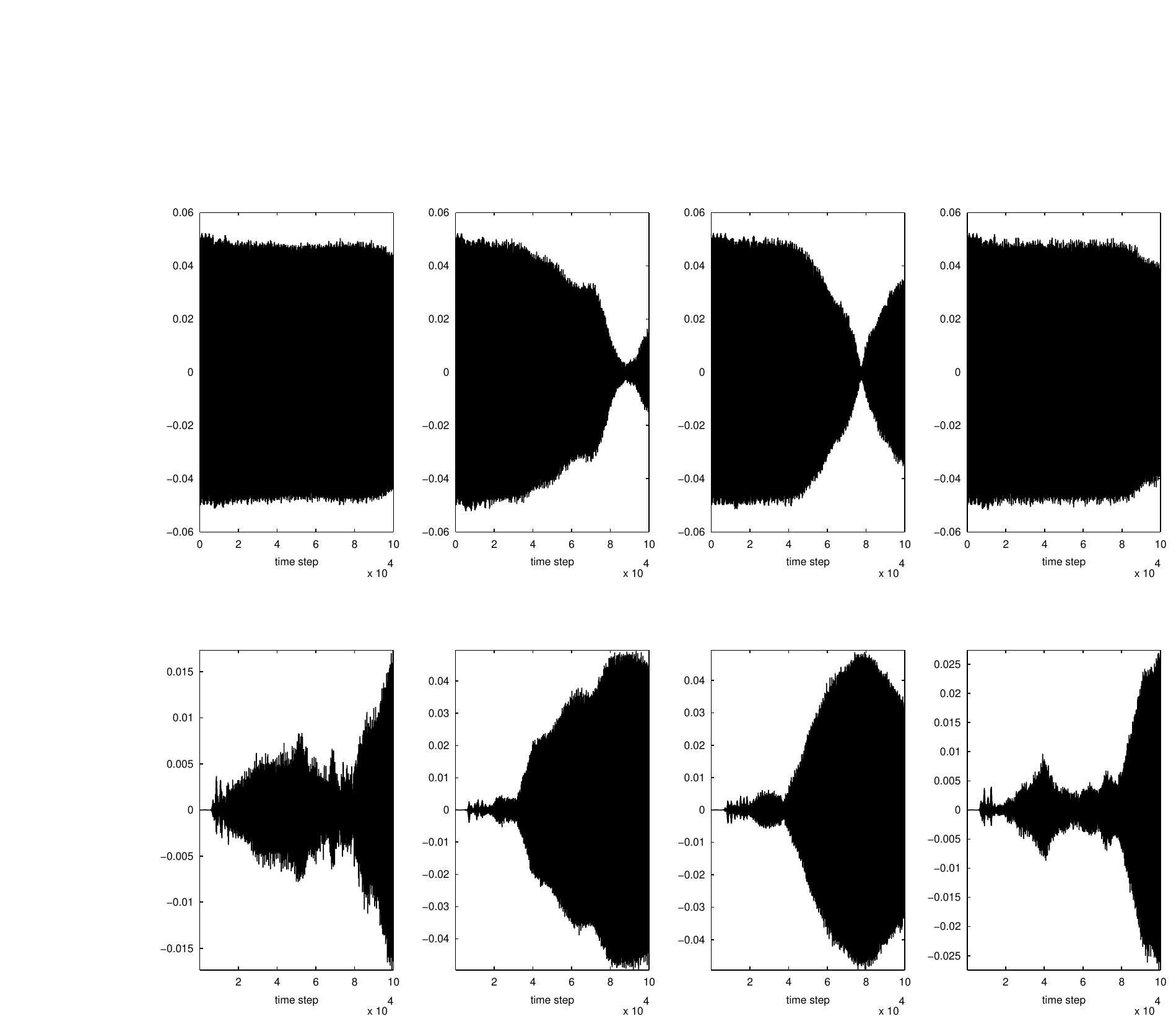}\multiput(-403,87)(101,0){4}{{\scriptsize $\eta_{(2)}$}}\multiput(-403,267)(101,0){4}{{\scriptsize $\eta_{(1)}$}}\put(-352,-20){(a)}\put(-251,-20){(b)}\put(-150,-20){(c)}\put(-50,-20){(d)}}\vspace{0.3in}
\caption{Planar instability. Transverse displacement at center of string (with $\alpha = 2\times 10^{-4}$), with fixed boundary conditions, plotted against time step, for various schemes for system ${\bf S}$: (a) scheme ${\bf s}^{(a)}$, (b) scheme ${\bf s}^{(b)}$, (c) scheme ${\bf s}^{(c)}$, and (d) scheme ${\bf s}^{(d)}$. The string is initialized using conditions \eqref{ex2}, with $\gamma_{1} = 0.05$ and $\gamma_{2} = 10^{-10}$; the time step is chosen as $h_{t}=1/10$. Both transverse polarizations are shown, (1) top row, and (2) bottom row. }  
\label{wifig}
\end{figure}

\begin{figure}[h!!!]
\centerline{\includegraphics[scale=0.85,clip,trim=0mm 100mm 23mm 40mm]{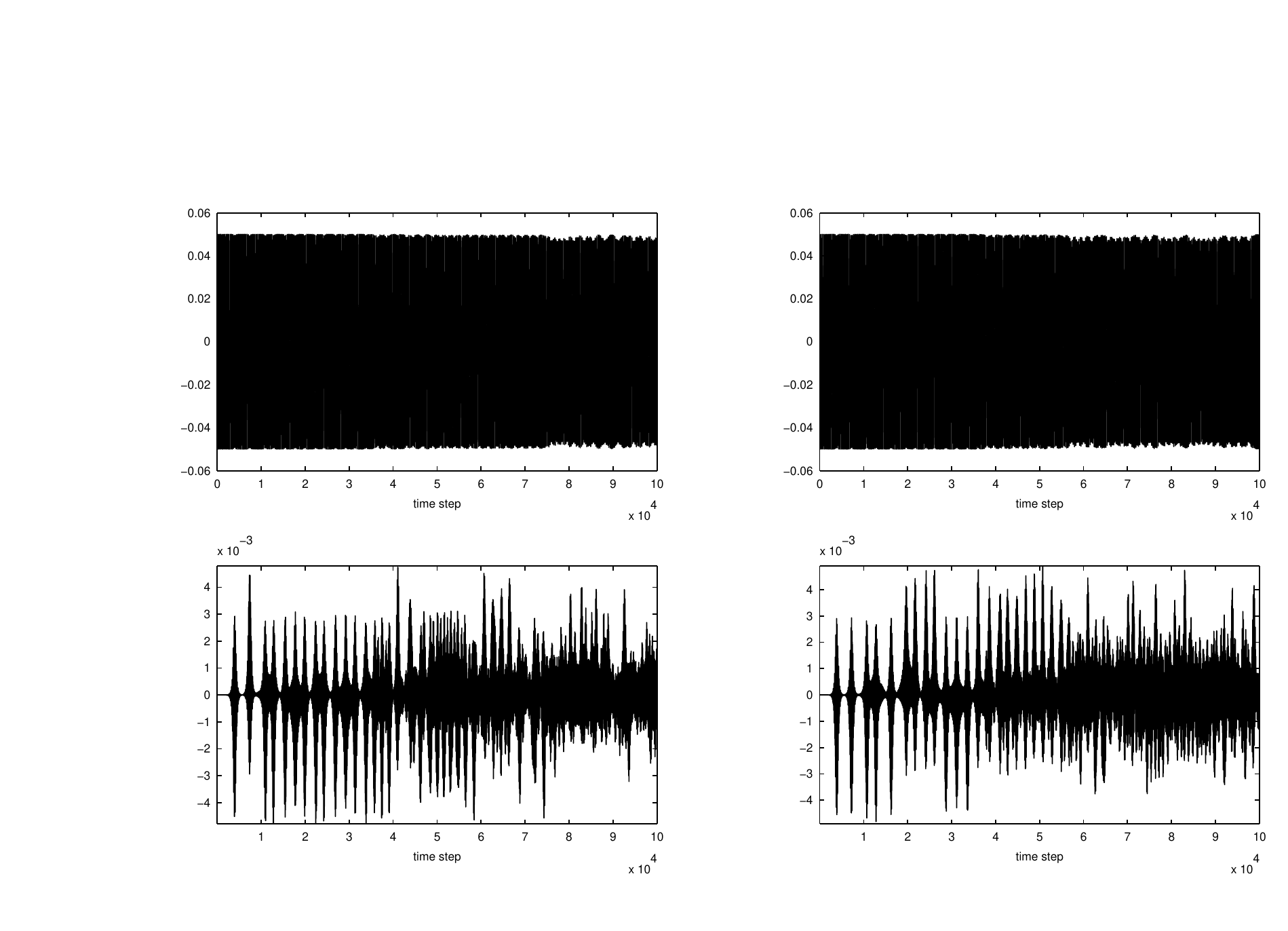}\multiput(-403,87)(217,0){2}{{\scriptsize $\eta_{(2)}$}}\multiput(-403,217)(217,0){2}{{\scriptsize $\eta_{(1)}$}}\put(-310,-20){(a)}\put(-90,-20){(b)}}\vspace{0.3in}
\caption{Planar instability. Transverse displacement at center of string (with $\alpha = 2\times 10^{-4}$), with fixed boundary conditions, plotted against time step, for schemes for system ${\bf K}$: (a) scheme ${\bf k}^{(a)}$, (b) and scheme ${\bf k}^{(b)}$. The string is initialized using conditions \eqref{ex2}, with $\gamma_{1} = 0.05$ and $\gamma_{2} = 10^{-10}$; the time step is chosen as $h_{t}=1/10$. Both transverse polarizations are shown, (1) top row, and (2) bottom row. }  
\label{wifig2}
\end{figure}

Consider a string under fixed conditions, with numerical initial conditions of the form
\begin{equation}
\label{ex2}
\eta_{(1),i}^0 = \eta_{(1),i}^1 = \gamma_{1}\sin(\pi ih_{x})\qquad \eta_{(2),i}^0 = \eta_{(2),i}^1 = \gamma_{2}\theta_{i}\sin(\pi ih_{x})
\end{equation}
where $\theta_{i}$ is a uniformly distributed random variable taking values over $(-1,1)$. In this case, the initial conditions correspond to a stationary modal distribution for $\eta_{(1)}$, of amplitude $\gamma_{1}$ and to a random perturbation for $\eta_{(2)}$, of amplitude $\gamma_{(2)}$. When $\gamma_{2}=0$, the state of the string remains in the $\eta_{(1)}$ polarization for all future time. But when $\gamma_{2}$ is not identically zero, even if very small, energy will be transferred from the $\eta_{(1)}$ polarization to the $\eta_{(2)}$ polarization, and whirling will occur. In Figs. \ref{wifig} and \ref{wifig2} are shown simulation results, again for the string of parameters given in the caption to Fig. \ref{fig1}, under the initial conditions \eqref{ex2}, using a very small value of $\gamma_{2}$ relative to $\gamma_{1}$ (values given in the captions). Both polarizations are shown in the figures. In the case of schemes ${\bf s}^{(\bullet)}$, energy is transferred nearly entirely from the $\eta_{(1)}$ polarization to the $\eta_{(2)}$, over approximately 100 000 time steps, and continues to oscillate henceforth. For schemes ${\bf k}^{(\bullet)}$, the transfer is much faster, and at the same time, limited, in the sense that the total energy transferred is much smaller. This planar instability, while physical, does not lead to explosive growth (indeed, angular momentum and energy conservation properties, where possessed by a given scheme, are not violated, as are stability conditions when implied by the latter form of conservation), but to a great deal of variation in the results of different schemes, particularly for schemes ${\bf s}^{(\bullet)}$. 

\subsection{Numerical Instability of Scheme ${\bf s}^{(b)}$}

Scheme ${\bf s}^{(b)}$ is perhaps the most interesting discussed here, in that it does indeed possess a discrete conserved energy, but this property does not lead to simple stability conditions. Instability which develops is inherently tied to ill-conditioning of the matrices${\bf A}_{{\bf s}^{(b)}}$ and ${\bf B}_{{\bf s}^{(b)}}$, from the matrix update form Eq. \eqref{updateform}, which leads to large numerical round off error. (Note in particular that of the matrices shown in Table \ref{tab5}, only ${\bf B}_{{\bf s}^{(b)}}$ is asymmetric.). 

\begin{figure}[h!!!]
\centerline{\includegraphics[scale=0.75,clip,trim=24mm 126mm 23mm 56mm]{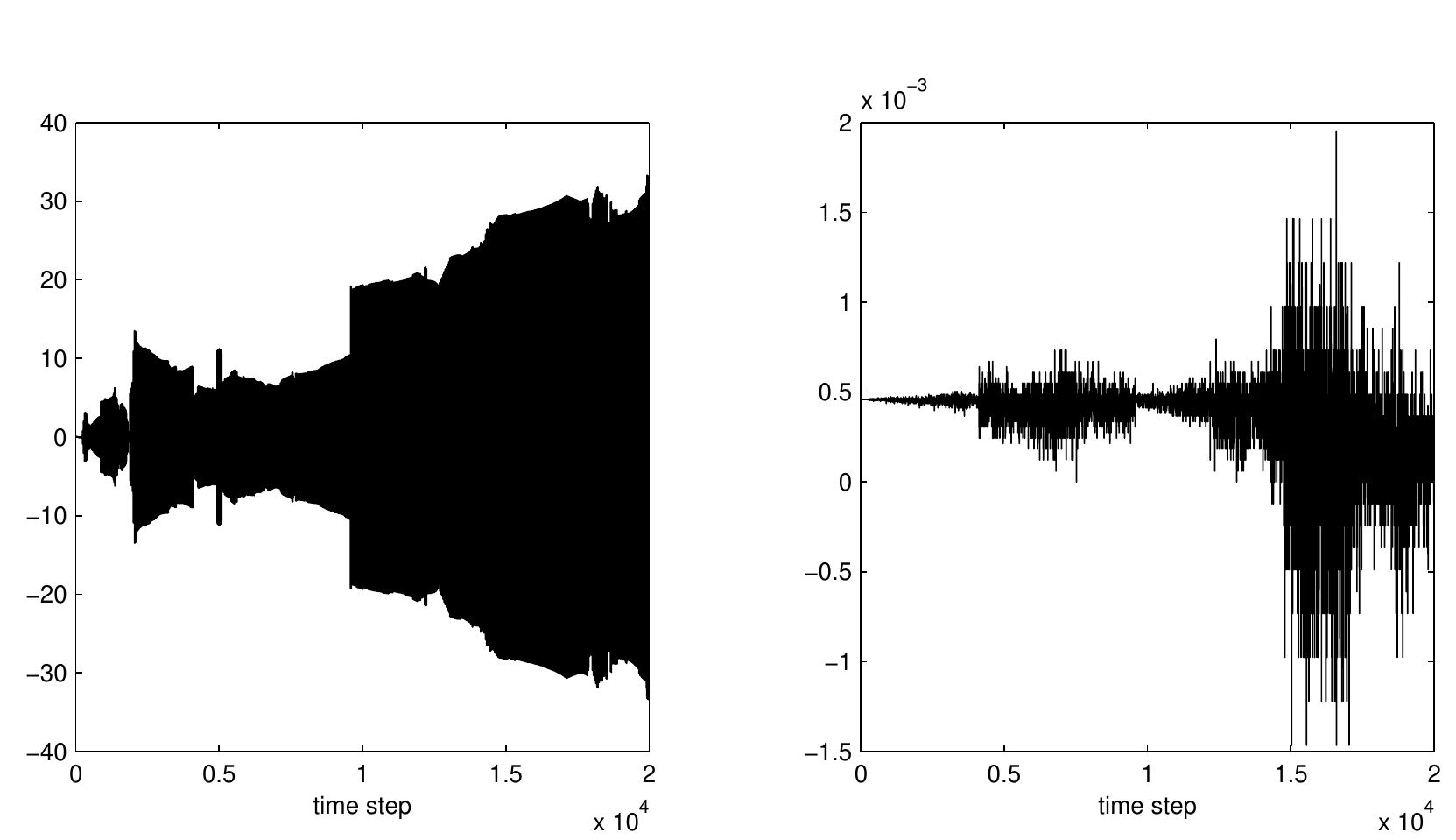}\put(-365,97){$\eta_{(1)}$}\put(-173,97){${\mathcal H}_{{\bf s}^{(b)}}$}\put(-277,-20){(a)}\put(-85,-20){(b)}}\vspace{0.3in}
\caption{Instability in system ${\bf s}^{(b)}$, applied to a string with $\alpha = 2\times 10^{-4}$, under initial condition \eqref{ex1}, with $\gamma_{1}=0.1$ and $\gamma_{2}=0$. The time step is chosen as $h_{t}=1/20$ s. (a), displacement $\eta_{(1)}$ and (b), energy ${\mathcal H}_{{\bf s}^{(b)}}$ plotted against time step.}  
\label{figu}
\end{figure}

Under low amplitude initial conditions, scheme ${\bf s}^{(b)}$ performs similarly to the other schemes for system ${\bf S}$, as discussed in Section \ref{numconssec}. Consider, though, planar initial conditions \eqref{ex1} with $\gamma_{1} = 0.1$ and $\gamma_{2}=0$. Though it is difficult to see in Figure \ref{figu}(b), the energy remains conserved and constant for several time steps, but then suddenly begins to fluctuate, and quickly begins taking on negative values. The displacement itself, shown in the left panel, experiences quantized jumps in amplitude, the largest coinciding with the moment at which the energy begins taking on negative values. If the simulation is allowed to proceed, similar discrete jumps in amplitude continue to occur; this type of instability is obviously of a very different nature from the exponential growth often seen in numerical schemes. Note in particular the interesting quantization effects in the fluctuations in energy.  

The schemes ${\bf s}^{(c)}$ and ${\bf s}^{(d)}$ are well-behaved under these conditions.

\section{Conclusions}

The main topic of this article has been the construction of numerical schemes for nonlinear strings which possess conservation properties; several different schemes for two distinct string models (systems ${\bf S}$ and ${\bf K}$) have been presented, which serve to highlight various distinctions among the schemes, which are summarized in Table \ref{tab10}. The angular momentum and energy conservation properties are generally independent; a scheme may possess one or the other, or both, or (though an example has not been presented here) neither. If energy conservation is implied by a scheme, stability conditions may follow, if the expression for energy can be shown to be generally positive for at least some range of choices of the material parameters and $\lambda$, the Courant number as defined by Eq. \eqref{Courantdef}. But, as in the case of scheme ${\bf s}^{(b)}$, this does not necessarily follow from energy conservation. Finally, there is the issue of computability, which is also independent of the other properties. Certain schemes are fully explicit (such as ${\bf s}^{(a)}$, ${\bf k}^{(a)}$ and ${\bf k}^{(b)}$), and others (namely ${\bf s}^{(b)}$, ${\bf s}^{(c)}$, ${\bf s}^{(d)}$ and ${\bf s}^{(e)}$) are formally implicit. There is a distinction to be made here, however, between schemes ${\bf s}^{(b)}$, ${\bf s}^{(c)}$, and ${\bf s}^{(d)}$, for which the implicit nature of the scheme is manifested as a linear system inversion performed at each time step, and ${\bf s}^{(e)}$, for which this is not possible. Thus, for this former set of schemes, there is at hand a simple proof of existence and uniqueness of solutions (i.e., the matrix update form Eq. \eqref{updateform}), whereas for the latter, there is not, and an implementation will require iterative methods for solution. It is important to reiterate that the simple implementation property of the former set of schemes is dependent on the form of the nonlinearity, which contains terms up to order three; for more general nonlinear forms (such as the exact nonlinear string system; see, e.g., \cite{Morse}), this is not the case. It is also tempting to conclude that, at least for a sufficiently complex nonlinear system such as ${\bf S}$, there is no fully explicit scheme which preserves energy; the same is clearly not true for the simpler system ${\bf K}$, for which the conservative scheme ${\bf k}^{(b)}$ may be written in an explicit form. 

\begin{table}
\caption{\label{tab10} Summary of properties of schemes ${\bf s}^{(\bullet)}$ and ${\bf k}^{(\bullet)}$. } 
\begin{center}\begin{scriptsize}
\begin{tabular}{l|l|l|l|l}
&\mbox{angular momentum conserving} & \mbox{energy conserving} & \mbox{stability conditions} & \mbox{existence/uniqueness}\\\hline\hline
${\bf s}^{(a)}$ & Yes& No & No & Yes\\\hline
${\bf s}^{(b)}$ & Yes& Yes& No & Yes\\\hline
${\bf s}^{(c)}$ & No & Yes& Yes& Yes\\\hline
${\bf s}^{(d)}$ & Yes& Yes& Yes& Yes\\\hline
${\bf s}^{(e)}$ & Yes& Yes& Yes& No\\\hline
${\bf k}^{(a)}$ & Yes& No & No & Yes\\\hline
${\bf k}^{(b)}$ & Yes& Yes& Yes& Yes\\
\end{tabular}
\end{scriptsize}\end{center}
\end{table}

Several interesting features of these methods were examined in numerical examples in Section \ref{numsec}. Planar instability of string motion was illustrated in Section \ref{planarsec}; under perturbed planar initial conditions, it is difficult to obtain consistent results across various different schemes for the same system; it may be possible to relate the resulting rate of oscillation to conservation of angular momentum, though this is a large topic, to be left to a future work. The largest issue, however, is that of discrete conservation itself, under finite-precision machine arithmetic. As was illustrated in Section \ref{numconssec}, energy and angular momentum may be conserved in a scheme exactly, in infinite precision, but when numerical round-off errors occur, conservation is lost. Under moderate vibration amplitudes, the resulting fluctuations will be small, often on the order of ``machine epsilon," but as the nonlinearity becomes stronger, these fluctuations can become large, even in floating point arithmetic, and numerical stability may be violated. Perhaps the best way of examining the effect of round-off will be through an examination of the matrix update form Eq. \eqref{updateform}, though this is a large separate topic which cannot be entered into in any detail here. The treatment of simulation techniques under finite-precision arithmetic is given scant attention in the literature, which is surprising; it is perhaps worth mentioning that this issue has been dealt with in the simulation of electrical networks (and also digital filtering), in particular by Fettweis (in the context of wave digital filters \cite{FettweisMain}) and Smith (for transmission lines and digital waveguides \cite{josbook}), essentially through the decomposition of conservative numerical methods into unitary matrix transformations. 

Finally, though the focus of this article has been on conservative methods for nonlinear string vibration, it should be clear that there are various more complex systems with nonlinearities of a similar simple form (i.e., resulting from a low-order series expansion); chief among these are the Berger and von Karman models of plate vibration \cite{Nayfeh}, \cite{Szilard}, which are, in a sense, direct analogues to systems ${\bf K}$ and ${\bf S}$ for the string, respectively, with added fourth-order terms modelling stiffness. The extension of the methods discussed here to these systems is immediate, and investigation is currently under way.

\appendix

\section{Linear Damping}
\label{losssec}
The introduction of linear damping terms to the systems under study here affects the analysis of difference schemes in only a very minor way. If such terms are included in system ${\bf S}$, Eqs. \eqref{Morse3} will be modified as
\begin{equation*}
\xi_{tt} = \phi_{x}-\sigma_{\xi}\xi_{t}\qquad\etab_{tt} = \psib_{x}-\sigma_{{\scriptsize\etab}}\etab_{t}
\end{equation*}
where $\sigma_{\xi}$ and $\sigma_{{\scriptsize\etab}}$ are non-negative constants, and where $\phi$ and $\psib$ are as defined in Eq. \eqref{Morse3}. The energetic analysis is similar to before, except that one will have 
\begin{equation*}
\frac{d}{dt}{\mathcal H}_{{\bf S}} = {\mathcal B}_{{\bf S}}-\sigma_{\xi}\|\xi_{t}\|^2-\sigma_{{\scriptsize\etab}}\|\etab_{t}\|^2
\end{equation*}
and if the boundary term ${\mathcal B}_{{\bf S}}$ vanishes, one has
\begin{equation*}
 \frac{d}{dt}{\mathcal H}_{{\bf S}}\leq 0\Longrightarrow{\mathcal H}_{{\bf S}}(t)\leq{\mathcal H}_{{\bf S}}(0)
\end{equation*}
Thus all bounds on solution size, discussed in Section \ref{boundsec}, remain unchanged. 

The conservation of angular momentum ${\mathcal A}_{{\bf S}}$, defined as before, is generalized to 
\begin{equation}
\frac{d}{dt}{\mathcal A}_{{\bf S}} = {\mathcal B}_{{\bf S},{\mathcal A}}-\sigma_{\scriptsize{\etab}}{\mathcal A}_{{\bf S}}
\end{equation}
In this case, under conservative boundary conditions ${\mathcal B}_{{\bf S},{\mathcal A}}=0$, one then has
\begin{equation}
{\mathcal A}_{{\bf S}}(t) = {\mathcal A}_{{\bf S}}(0)e^{-\sigma_{{\scriptsize\etab}}t}
\end{equation}
In other words, the angular momentum decays exponentially. 

Consider one of the difference schemes ${\bf s}^{(c)}$, ${\bf s}^{(d)}$ or ${\bf s}^{(e)}$, which is both energy conserving and for which CFL-like conditions exist for the positivity of the discrete energy. For any of these schemes, an extension to the lossy case may be discretized as
\begin{equation*}
\delta_{t+}\delta_{t-}\xi = \delta_{x+}\phi_{{\bf s}^{(\bullet)}}-\sigma_{\xi}\delta_{to}\xi\qquad\delta_{t+}\delta_{t-}\etab = \delta_{x+}\psib_{{\bf s}^{(\bullet)}}-\sigma_{{\scriptsize\etab}}\delta_{to}\etab
\end{equation*}
Discrete energy conservation generalizes simply to 
\begin{equation*}
\delta_{t+}{\mathcal H}_{{\bf s}^{(\bullet)}} = {\mathcal B}_{{\mathcal H}, {\bf s}^{(\bullet)}}-\sigma_{\xi}\|\delta_{to}\xi\|_{{\mathcal D}}^2-\sigma_{{\scriptsize\etab}}\|\delta_{to}\etab\|_{{\mathcal D}}^2
\end{equation*}
and if the boundary terms vanish, one then has
\begin{equation*}
 \delta_{t+}{\mathcal H}_{{\bf s}^{(\bullet)}}\leq 0\Longrightarrow{\mathcal H}_{{\bf s}^{(\bullet)}}^{n}\leq{\mathcal H}_{{\bf s}^{(\bullet)}}^{0} 
\end{equation*}
and all derived bounds on solution size again hold as before. 

For schemes ${\bf s}^{(a)}$ and ${\bf s}^{(d)}$, it is simple to show that
\begin{equation}
\delta_{t+}{\mathcal A}_{{\bf s}^{\bullet}} = {\mathcal B}_{{\mathcal A}, {\bf s}^{\bullet}}-\sigma_{{\scriptsize\etab}}\mu_{t+}{\mathcal A}_{{\bf s}^{\bullet}}
\end{equation}
and if the boundary terms vanish, one has a geometric decay in the angular momentum of the form
\begin{equation}
{\mathcal A}_{{\bf s}^{\bullet}}^{n} = \frac{1-\frac{h_{t}\sigma_{{\scriptsize\etab}}}{2}}{1+\frac{h_{t}\sigma_{{\scriptsize\etab}}}{2}}{\mathcal A}_{{\bf s}^{\bullet}}^{n-1}
\end{equation}
(Notice, however, that if 
\begin{equation}
h_{t}\geq 2/\sigma_{{\scriptsize\etab}}
\end{equation}
the decay will be highly oscillatory, and unphysical.) For schemes ${\bf s}^{(b)}$ and ${\bf s}^{(e)}$, such exact geometric decay does not follow. 

When a linear damping term is added to system ${\bf K}$, it is easy to show that a similar monotonic decrease in total energy, as well as an exponential decay in angular momentum also follows. For both schemes ${\bf k}^{(a)}$ and ${\bf k}^{(b)}$, when the linear damping term is added, similarly to for schemes ${\bf s}^{(\bullet)}$ above, angular momentum will decrease geometrically, and for scheme ${\bf k}^{(b)}$, discrete energy will decrease monotonically. 

\section{Improved Stability Conditions for Schemes ${\bf s}^{(\bullet)}$}
\label{impsec}
For any of the schemes for system ${\bf S}$ which are conservative and for which a global stability condition may be derived (namely schemes ${\bf s}^{(c)}$, ${\bf s}^{(d)}$, and ${\bf s}^{(e)}$), it is true that the condition on the time step $h_{t}$, for a given grid spacing $h_{x}$ can be exceedingly small, from conditions such as \eqref{CFLcond} and \eqref{CFLcond3}. On the other hand, these schemes are already implicit, and there is thus no loss in efficiency in generalizing the treatment of the linear part of system ${\bf S}$ to improve on these strict bounds on $h_{t}$.  

Turning attention to Table \ref{tab1}, any of schemes ${\bf s}^{(\bullet)}$ may be generalized directly as
\begin{subequations}
\label{impsys}
\begin{eqnarray}
\delta_{t+}\delta_{t-}\xi &=& \delta_{x+}\phi_{{\bf s}^{(\bullet)},\tau}\\
\delta_{t+}\delta_{t-}\etab &=& \delta_{x+}\psib_{{\bf s}^{(\bullet)},\nu}
\end{eqnarray}
\end{subequations}
where $\phi_{{\bf s}^{(\bullet)},\tau}$ and $\psib_{{\bf s}^{(\bullet)},\nu}$ are defined in terms of the free parameters $\tau$ and $\nu$ as
\begin{equation*}
\phi_{{\bf s}^{(\bullet)},\tau} = \phi_{{\bf s}^{(\bullet)}}+\tau\left(p-\mu_{t+}\mu_{t-}p\right)\qquad \psib_{{\bf s}^{(\bullet)},\nu} = \psib_{{\bf s}^{(\bullet)}}+\alpha\nu\left(\etab-\mu_{t+}\mu_{t-}\etab\right)
\end{equation*}
Schemes \eqref{impsys} reduce to the forms shown in Table \ref{tab1} when $\tau=\nu=0$. 

Consider first the energetic analysis of the generalized systems \eqref{impsys}; this is essentially the same as that carried out in Section \ref{egyconssec}, except for the new linear terms parameterized by $\tau$ and $\nu$. For those schemes with an energy conservation property, namely  ${\bf s}^{(b)}$, ${\bf s}^{(c)}$, ${\bf s}^{(d)}$ and ${\bf s}^{(e)}$, the potential energy under this generalization is modified to
\begin{equation*}
{\mathcal V}_{{\bf s}^{(\bullet)},\tau,\nu} = {\mathcal V}_{{\bf s}^{(\bullet)}}+\frac{\tau h_{t}^{2}}{8}\|\delta_{t-}p\|_{{\mathcal D}^{+}}^{2}+\frac{\alpha\nu h_{t}^{2}}{8}\|\delta_{t-}q\|_{{\mathcal D}^{+}}^{2}
\end{equation*}
In other words, for positive $\tau$ and $\nu$, the additional contribution to the potential energy is positive. 

Consider now, as an example, the effect on scheme ${\bf s}^{(d)}$. The potential energy is generalized from the expression \eqref{vd} to
\begin{eqnarray*}
{\mathcal V}_{{\bf s}^{(d)},\tau,\nu} &=& \frac{\alpha}{2}\left(\|\mu_{t-}p\|^2_{{\mathcal D}^{+}}+ \|\mu_{t-}{\bf q}\|^2_{{\mathcal D}^{+}}\right)+\frac{1-\alpha}{2}\left(\|\mu_{t-}p+\frac{1}{2}{\bf q}^{T}e_{t-}{\bf q}\|^{2}_{{\mathcal D}^{+}}\right)\\
\qquad&& +\frac{(\tau-1)h_{t}^2}{8}\|\delta_{t-}p\|^2_{{\mathcal D}^{+}}+\frac{\alpha(\nu-1)h_{t}^2}{8}\|\delta_{t-}{\bf q}\|^2_{{\mathcal D}^{+}}\\
\end{eqnarray*}
The expression for the kinetic energy ${\mathcal T}_{{\bf s}^{(\bullet)}}$ as given in Table \ref{tab3} remains unchanged. 
Clearly, if $\tau\geq 1$ and $\nu\geq 1$, then the potential energy is always positive, as will be the energy ${\mathcal H}_{{\bf s}^{(d)}, \tau,\nu}$. Under these conditions, the generalized algorithm is stable for any choice of time step $h_{t}$. If $\tau\leq 1$, or $\nu\leq 1$, the analysis is very similar to that carried out previously. One may obtain similar improvements on the bound on the time step for schemes ${\bf s}^{(d)}$ and ${\bf s}^{(e)}$, as well as ${\bf k}^{(a)}$ and ${\bf k}^{(b)}$; scheme ${\bf s}^{(b)}$ can still be shown to have a conserved energy which is non-positive for at least some choices of state variables. 

As mentioned earlier in this section, the generalized schemes will not require more operations per time step; some modifications to the update matrices given in Section \ref{implementationsec} will be necessary, but the sparsity remains nearly the same.

\section{Energy-conserving Spectral Method for System {\bf K}}
\label{specsec}
Due to the special form of the nonlinearity in system ${\bf K}$, an alternative analysis is possible using  spatial Fourier series expansion techniques \cite{Carrier}, \cite{Dickey}, \cite{Dickey80}; such analysis leads naturally to the construction of highly accurate spectral-type numerical solution methods \cite{trefethenspec}, \cite{Fornberg}, which like the simpler difference schemes discussed in the main body of this article, are conservative. The same is not true for system ${\bf S}$. 

Consider system ${\bf K}$ under fixed boundary conditions. An expansion for \etab\, of the form 
\begin{equation*}
\etab (x,t) = \sum_{k=1}^{\infty}\hat{\etab}_{k}(t)\sin(\pi k x)
\end{equation*}
where the time-dependent vector expansion coefficients are given by $\hat{\etab}_{k}(t)$, for $k=1,\hdots,\infty$, thus satisfies the boundary conditions automatically. One may then rewrite system ${\bf K}$ as the infinite system of ordinary differential equations 
\begin{equation*}
\frac{d^2\hat{\etab}_{k}}{dt^2} = -\alpha{\mathcal G}\pi^2 k^2\hat{\etab}_{k}\qquad k = 1,\hdots,\infty
\end{equation*}
This may be time-discretized immediately, and the infinite system of equations truncated to $M$ terms to form system ${\bf k}^{(s)}$, defined as
\begin{equation}
\label{kspec}
\delta_{t+}\delta_{t-}\hat{\etab}_{k} = -\alpha{\mathcal G}_{{\bf k}^{(s)}}\pi^2 k^2\hat{\etab}_{k}\qquad k = 1,\hdots,M
\end{equation}
where the form of ${\mathcal G}_{{\bf k}^{(s)}}$ under discretization is left unspecified for the moment. Notice that the approximation to the spatial derivative operators above is spectrally accurate, and is exact in the limit as $M$ becomes large.

Introducing the inner product of two sets of vector expansion coefficients $\hat{\bf f}_{k}$ and $\hat{\bf g}_{k}$ of dimension $M$ by 
\begin{equation*}
\langle \hat{{\bf f}},\hat{{\bf g}}\rangle_{[1,M]} = 2\sum_{k=1}^{M}\hat{{\bf f}}_{k}^{T}\hat{{\bf g}}_{k}
\end{equation*}
and the associated norm by
\begin{equation*}
\|\hat{{\bf f}}\|_{[1,M]} = \langle \hat{{\bf f}},\hat{{\bf f}}\rangle_{[1,M]}^{1/2}
\end{equation*}
one may then take the inner product of Eq. \eqref{kspec} with $\delta_{to}\hat{\etab}$ to get
\begin{equation*}
\langle \delta_{to}\hat{\etab},\delta_{t+}\delta_{t-}\hat{\etab}\rangle_{[1,M]} = -\alpha{\mathcal G}_{{\bf k}^{(s)}}\langle \delta_{to}\hat{\etab},\pi^2 k^2\hat{\etab}\rangle_{[1,M]}
\end{equation*}
or
\begin{equation*}
\delta_{t+}\left(\frac{1}{2}\| \delta_{t-}\hat{\etab}\|_{[1,M]}^2\right) +\frac{\alpha{\mathcal G}_{{\bf k}^{(s)}}}{2}\delta_{t+}\langle\pi k\hat{\etab},\pi k e_{t-}\hat{\etab}\rangle_{[1,M]}=0
\end{equation*}
Now, in analogy with scheme ${\bf k}^{(b)}$, one may define ${\mathcal G}_{{\bf k}^{(s)}}$ as
\begin{equation*}
{\mathcal G}_{{\bf k}^{(s)}} = 1+\frac{1}{2\alpha}\mu_{t+}\langle\pi k\hat{\etab},\pi k e_{t-}\hat{\etab}\rangle_{[1,M]}
\end{equation*}
and one again arrives at an expression for conserved energy as
\begin{equation*}
\delta_{t+}{\mathcal H}_{{\bf k}^{(s)}} = \delta_{t+}\left({\mathcal T}_{{\bf k}^{(s)}}+{\mathcal V}_{{\bf k}^{(s)}}\right)=0
\end{equation*}
with
\begin{eqnarray*}
{\mathcal T}_{{\bf k}^{(s)}} &=& \frac{1}{2}\| \delta_{t-}\hat{\etab}\|_{[1,M]}^2\\
{\mathcal V}_{{\bf k}^{(s)}} &=& \frac{\alpha}{2}\langle\pi k\hat{\etab},\pi k e_{t-}\hat{\etab}\rangle_{[1,M]}\left(1+\frac{1}{4\alpha}\langle\pi k\hat{\etab},\pi k e_{t-}\hat{\etab}\rangle_{[1,M]}\right)
\end{eqnarray*}

Conditions for positivity may be arrived at by rewriting the expression for potential energy as
\begin{eqnarray*}
{\mathcal V}_{{\bf k}^{(s)}} &=& \frac{\alpha\pi^2}{2}\left(\|k\mu_{t-}\hat{\etab}\|_{[1,M]}^2-\frac{h_{t}^2}{4}\|k\delta_{t-}\hat{\etab}\|_{[1,M]}^2\right)+\frac{1}{8}\langle\pi k\hat{\etab},\pi k e_{t-}\hat{\etab}\rangle_{[1,M]}^2\notag\\
&\geq& \frac{\alpha\pi^2}{2}\left(\|k\mu_{t-}\hat{\etab}\|_{[1,M]}^2-\frac{h_{t}^2 M^2}{4}\|\delta_{t-}\hat{\etab}\|_{[1,M]}^2\right)+\frac{1}{8}\langle\pi k\hat{\etab},\pi k e_{t-}\hat{\etab}\rangle_{[1,M]}^2
\end{eqnarray*}
which gives a lower bound for the total energy as
\begin{eqnarray*}
{\mathcal H}_{{\bf k}^{(s)}}&\geq& \frac{\alpha\pi^2}{2}\|k\mu_{t-}\hat{\etab}\|_{[1,M]}^2+\frac{1}{8}\langle\pi k\hat{\etab},\pi k e_{t-}\hat{\etab}\rangle_{[1,M]}^2\notag\\
&+& \left(\frac{1}{2}-\frac{\alpha h_{t}^2\pi^2 M^2}{8}\right)\|\delta_{t-}\hat{\etab}\|_{[1,M]}^2
\end{eqnarray*}
and the positivity condition is easily read off as
\begin{equation*}
h_{t}\leq \frac{2}{\pi M}\sqrt{\frac{1}{\alpha}}
\end{equation*}

Given this positivity condition, bounds on the solution size (i.e., bounds on the norm of $\hat{\etab}$) may be derived exactly as in Section \ref{boundsec}. Bounds on $\hat{\etab}$ may be simply related to bounds on \etab\,
 itself through an application of Parseval's Theorem \cite{Horn}. Conservation of angular momentum also holds for system ${\bf k}^{(s)}$, and is trivial to show, through an inner product of system \eqref{kspec} with $\tilde{\hat{\etab}}$.


\bibliographystyle{unsrt}
\bibliography{kirchoffbib}

\end{document}